\documentclass[reqno]{amsart}

\usepackage{amssymb,amsmath,amsthm,amsfonts}
\usepackage{graphicx}

\theoremstyle{plain}
\newtheorem{theorem}{Theorem}[section]
\newtheorem*{theorem*}{Theorem}
\newtheorem{proposition}[theorem]{Proposition}
\newtheorem*{proposition*}{Proposition}
\newtheorem{corollary}[theorem]{Corollary}
\newtheorem*{corollary*}{Corollary}
\newtheorem{lemma}[theorem]{Lemma}
\newtheorem*{lemma*}{Lemma}

\theoremstyle{definition}
\newtheorem{definition}[theorem]{Definition}
\newtheorem*{definition*}{Definition}

\newtheorem*{definitions*}{Definitions}

\newtheorem*{example*}{Example}

\newtheorem*{examples*}{Examples}

\newtheorem*{exercise*}{Exercise}

\newtheorem*{exercises*}{Exercises}

\theoremstyle{remark}

\newtheorem*{remark*}{Remark}

\newtheorem*{remarks*}{Remarks}

\newtheorem*{note*}{Note}

\newtheorem*{notes*}{Notes}

\newtheorem*{notation*}{Notation}

\numberwithin{equation}{section}

\newcommand{\R}{\mathbb{R}} 
\newcommand{\N}{\mathbb{N}} 
\newcommand{\Om}{\Omega} \newcommand{\w}{\omega}
\newcommand{\eps}{\varepsilon} \newcommand{\di}{\textsf{d}}
\newcommand{\wt}{\widetilde} \newcommand{\n}[1]{\| #1 \|}
\newcommand{\wh}{\widehat} \DeclareMathOperator{\car}{card}
\DeclareMathOperator{\cls}{cls}

\begin{document}
\title[Exponential ordering for NFDEs with non-autonomous $D$-operator]
{Exponential ordering for neutral functional differential equations with
  non-autonomous linear $\mathbf D$-operator} \author[R. Obaya]{Rafael
  Obaya} \author[V. M. Villarragut]{V\'{\i}ctor M. Villarragut}
\address{Departamento de Matem\'{a}tica Aplicada, Escuela de
  Ingenier{\'\i}as Industriales, Universidad de Valladolid, 47011
  Valladolid, Spain.}  \email[Rafael Obaya]{rafoba@wmatem.eis.uva.es}
\email[V\'{\i}ctor M. Villarragut]{vicmun@wmatem.eis.uva.es} \thanks{The
  authors were partly supported by Junta de Castilla y Le\'{o}n under
  project VA060A09, and Ministerio de Ciencia e Innovaci{\'o}n under
  project MTM2008-00700/MTM.}  \date{}
\dedicatory{Dedicated to Professor Russell A. Johnson on the occasion of
his 60th birthday}
\begin{abstract}
  We study neutral functional differential equations with stable linear
  non-autonomous $D$-operator. The operator of convolution $\wh D$
  transforms $BU$ into $BU$. We show that, if $D$ is stable, then $\wh D$
  is invertible and, besides, $\wh D$ and $\wh D^{-1}$ are uniformly
  continuous for the compact-open topology on bounded sets. We introduce a
  new transformed exponential order and, under convenient assumptions, we
  deduce the 1-covering property of minimal sets. These conclusions are
  applied to describe the amount of material in a class of compartmental
  systems extensively studied in the literature.
\end{abstract}
\keywords{Non-autonomous dynamical systems, monotone skew-product
  semiflows, neutral functional differential equations, infinite delay,
  compartmental systems} \subjclass{37B55, 34K40, 34K14}
\renewcommand{\subjclassname}{\textup{2000} Mathematics Subject
  Classification}

\vspace*{-1.5cm}

\noindent\begin{tabular}{|l|}
  \hline
  Obaya, R., Villarragut, V.M.,\\
  Exponential Ordering for Neutral Functional Differential Equations\\
With Non-Autonomous Linear D-Operator. J Dyn Diff Equat 23, 695–725 (2011).\\
https://doi.org/10.1007/s10884-011-9210-9\\
  \copyright ~Springer\\
  \hline
\end{tabular}

\vspace{1cm}

\maketitle
\section{Introduction}
The present work is one of a series of papers devoted to the study of the
long term behavior of the solutions of non-autonomous monotone differential
equations and to provide a local or global qualitative description of the
corresponding phase spaces (see Shen and Yi~\cite{shyi}, Jiang and
Zhao~\cite{jizh}, and Novo, Obaya and Sanz~\cite{noos} among others). This
paper continues the theory initiated by Mu{\~n}oz-Villarragut, Novo and
Obaya~\cite{muno}, and Novo, Obaya and Villarragut~\cite{noov}, where a
dynamical theory for non-autonomous functional differential equations (FDEs
for short) with infinite delay and for non-autonomous neutral functional
differential equations (NFDEs for short) with infinite delay and stable
autonomous linear $D$-operator was developed. We now consider
non-autonomous NFDEs with non-autonomous linear $D$-operator. In this
situation, the main conclusions of the previous papers do not remain valid
and, thus, the extension of the theory requires to solve new and
challenging problems, which becomes the main objective of this paper. We
also introduce an alternative definition of the exponential ordering which
can be applied in the present context with more generality preserving the
dynamical behavior of the previous theory.\par
We assume some recurrence properties on the temporal variation of the
NFDEs; thus, its solutions induce a skew-product semiflow with a minimal
flow $(\Om,\sigma,\R)$ on the base. In particular, the uniform almost
periodic and almost automorphic cases are included in this formulation. The
skew-product formalism permits the study of the trajectories using methods
of ergodic theory and topological dynamics. We introduce the phase space
$BU\subset C((-\infty,0],\R^m)$ of bounded and uniformly continuous
functions with the supremum norm, where the standard theory provides
existence, uniqueness and continuous dependence of the solutions. For each
$r>0$, we denote $B_r=\{x\in BU:\n{x}_\infty\leq r\}$. In our present
setting, every bounded trajectory is relatively compact for the
compact-open metric topology, the restriction of the semiflow to its
omega-limit sets is continuous for the metric and it admits a flow
extension. In fact, the main conclusions of this theory are obtained by the
application of the metric topology on compact invariant sets. We introduce
a non-autonomous operator $D:\Om\times BU\to\R^m$ which is linear and
continuous for the norm in its second variable, continuous for the metric
topology on the bounded subsets of $BU$, and atomic at zero. We obtain an
integral representation $D(\w,x)=\int_{-\infty}^0[d\mu(\w)]x$,
$(\w,x)\in\Om\times BU$, where $\mu(\w)$ is an $m\times m$ matrix of real
Borel regular measures with finite total variation and
$\det(\mu(\w)(\{0\}))\neq 0$ for every $\w\in\Om$. We can define $\wh
D:\Om\times BU\to\Om\times BU$ by $\wh D(\w,x)=(\w,\wh D_2(\w,x))$ with
$\wh D_2(\w,x)(t)=\int_{-\infty}^0[d\mu(\w{\cdot}t)(s)]x(t+s)$ for all
$t\leq 0$. Then $\wh D_2$ is well-defined, it is linear and continuous for
the norm in its second variable for all $\w\in\Om$ and it is uniformly
continuous on $\Om\times B_r$ for every $r>0$ when we consider the
restriction of the compact-open topology in the second factor. We prove
that, if $D$ is stable in the sense of Hale~\cite{hale}, and Hale and
Verduyn-Lunel~\cite{have}, then $\wh D$ is invertible, $(\wh
D^{-1})_2(\w,{\cdot})$ is linear and continuous on $BU$ and $\wh D^{-1}$ is
uniformly continuous on $\Om\times B_r$ for every $r>0$ when we take the
restriction of the compact-open topology in the second factor. This
behavior characterizes the stability of the operator. In fact $D$ is stable
if and only if $\wh D$ is invertible and $\wh D^{-1}$ is continuous for the
product metric topology on bounded subsets. In addition, if $D$ is stable
and the flow $(\Om,\sigma,\R)$ is almost periodic, then $\wh D$ and $\wh
D^{-1}$ are uniformly continuous on $\Om\times B_r$ for all $r>0$ when we
take the norm in the second factor. In this almost periodic context, our
theory is completely analogous to the one obtained in~\cite{muno} for the
autonomous operators $D$ and $\wh D$.\par
The exponential order was extensively studied by Smith and
Thieme~\cite{smth}, \cite{smth2} for autonomous functional differential
equations with finite delay. Assuming that the vector field satisfies a
strong version of a non-standard monotonicity condition and a supplementary
irreducibility assumption, they proved that the induced semiflow is
eventually strongly monotone in the phase space of Lipschitz functions and
strongly order preserving in the phase space of continuous functions, from
which they deduced the quasiconvergence of the trajectories. Krisztin and
Wu~\cite{krwu} studied the dynamical properties of scalar neutral
functional differential equations with finite delay which induce a monotone
semiflow for the exponential ordering, and Wu and Zhao~\cite{wuzh} analyzed
the same question in a class of evolutionary equations with applications to
reaction-diffusion models. Skew-product semiflows generated by families of
non-autonomous FDEs with infinite delay and NFDEs with infinite delay and
autonomous stable $D$-operator which are monotone for the exponential order
were studied in~\cite{noov}. Assuming some adequate hypotheses including
the boundedness, relative compactness and uniform stability for the order
on bounded subsets of the trajectories, they deduced that the omega-limit
set of either each trajectory or at least each trajectory with Lipschitz
continuous initial datum is a copy of the base, i.e. it reproduces exactly
the dynamics exhibited by the time variation of the equation. It is
important to point out that no irreducibility conditions are required, and
hence the conclusion can be applied to general problems under natural
physical conditions.\par
In the setting considered in this paper, it is obvious that the
invertibility of the operator $\wh D$ allows us to transform the original
NFDEs with a time dependent $D$-operator into FDEs to which the conclusions
in~\cite{noov} can be applied. Thus we can define a \emph{transformed
  exponential order} on our original phase space: we say that a semiflow is
monotone for the transformed exponential order when the transformed
semiflow is monotone for the exponential order introduced
in~\cite{noov}. This transformed exponential order has the same dynamical
properties as the direct exponential order; under natural conditions of
relative compactness and stability, the omega-limit sets are copies of the
base. Two different approaches may be taken; on the one hand, we can use
the exponential order for the transformed FDE on a compact interval of
$(-\infty,0]$ and the standard ordering on its complementary set; on the
other hand, we can consider the exponential order for the transformed FDE
on the complete interval $(-\infty,0]$. Each approach yields different
dynamical consequences.\par
The foregoing conclusions are applied to the study the evolution of the
amount of material within the compartments of some non-autonomous neutral
compartmental systems, extending previous results of the
literature. Compartmental systems have been used as mathematical models for
the study of the dynamical behavior in biological and physical sciences
which depend on local mass balance conditions (see Jacquez~\cite{jacq},
Jacquez and Simon~\cite{jasi}, \cite{jasi2}, Haddad, Chellaboina and
Hiu~\cite{hach}, and the references therein). The papers \cite{krwu},
\cite{muno}, \cite{noov}, Arino and Bourad~\cite{arbo}, and Wu and
Freedman~\cite{wufr} apply dynamical methods in order to study monotone
neutral compartmental systems. In this work, we consider a general class of
compartmental systems and solve, under convenient hypotheses, the same kind
of problems by using the transformed exponential order. Assuming the
monotonicity and the existence of an adequate bounded semitrajectory, we
conclude that either every trajectory or at least every trajectory with
Lipschitz continuous transformed initial datum is bounded and uniformly
stable for the transformed order on bounded sets.\par
We analyze in detail some specific closed neutral compartmental systems
whose explicit expression ensures the stability of $D$. Now the total mass
of the system is invariant and we obtain precise conditions guaranteeing
the monotonicity for the transformed exponential order of the
semiflow. Under these conditions, the omega-limit set of either every
trajectory or at least every trajectory with Lipschitz continuous
transformed initial datum is a copy of the base. We also give some
alternative conditions leading to the monotonicity of the semiflow for the
direct exponential order, which requires the differentiability of the
coefficients of $D$ rather than just their continuity, used for the
transformed exponential order. Thus its applicability becomes much more
restrictive. In addition the transformed exponential order is also more
advantageous when dealing with rapidly oscillating coefficients of $D$; as
a consequence, this is the natural exponential order when dealing with
NFDEs with recurrent linear $D$-operator.\par
We now briefly describe the structure of the paper. Some basic notions on
topological dynamics used throughout the rest of the paper are stated in
Section~\ref{section_preliminaries}. In Section~\ref{section_operator}, we
study the invertibility of $\wh D$ and show that, when $D$ is stable, the
theories of $\wh D$ and $\wh D^{-1}$ are symmetric. In
Section~\ref{section_copy_base}, we introduce the transformed exponential
order and prove that the main conclusions of Sections~4 and~5
in~\cite{noov} remain valid for this order. Section~\ref{section_comp_sys}
studies the boundedness and stability for the order of the trajectories of
a quite general class of non-autonomous compartmental systems. Finally, in
Section~\ref{section_examples} we analyze some specific closed neutral
compartmental systems, we characterize the relatively compact trajectories
and prove that the minimal sets are copies of the base.
\section{Some preliminaries}\label{section_preliminaries}
Let $(\Om,d)$ be a compact metric space. A real {\em continuous flow \/}
$(\Om,\sigma,\R)$ is defined by a continuous mapping $\sigma: \R\times \Om
\to \Om,\; (t,\w)\mapsto \sigma(t,\w)$ satisfying
\begin{enumerate}
  \renewcommand{\labelenumi}{(\roman{enumi})}
\item $\sigma_0=\text{Id},$
\item $\sigma_{t+s}=\sigma_t\circ\sigma_s$ for each $s$, $t\in\R$,
\end{enumerate}
where $\sigma_t(\w)=\sigma(t,\w)$ for all $\w \in \Om$ and $t\in \R$. The
set $\{ \sigma_t(\w) : t\in\R\}$ is called the {\em orbit\/} or the {\em
  trajectory\/} of the point $\w$. We say that a subset $\Om_1\subset \Om$
is {\em $\sigma$-invariant\/} if $\sigma_t(\Om_1)=\Om_1$ for every
$t\in\R$.  A subset $\Om_1\subset \Om$ is called {\em minimal \/} if it is
compact, $\sigma$-invariant and its only nonempty compact
$\sigma$-invariant subset is itself. Every compact and $\sigma$-invariant
set contains a minimal subset; in particular it is easy to prove that a
compact $\sigma$-invariant subset is minimal if and only if every
trajectory is dense. We say that the continuous flow $(\Om,\sigma,\R)$ is
{\em recurrent\/} or {\em minimal\/} if $\Om$ is minimal.\par
The flow $(\Om,\sigma,\R)$ is {\em distal\/} if for any two distinct points
$\w_1,\,\w_2\in\Om$ the orbits keep at a positive distance, that is,
$\inf_{t\in \R}d(\sigma(t,\w_1),\sigma(t,\w_2))>0$. The flow
$(\Om,\sigma,\R)$ is {\em almost periodic\/} when for every $\varepsilon >
0 $ there is a $\delta >0$ such that, if $\w_1$, $\w_2\in\Om$ with
$d(\w_1,\w_2)<\delta$, then $d(\sigma(t,\w_1),\sigma(t,\w_2))<\varepsilon$
for every $t\in \R$. Equivalently, the flow $(\Om,\sigma,\R)$ is almost
periodic if the family $\{\sigma_t\}_{t\in\R}$ is equicontinuous. If
$(\Om,\sigma,\R)$ is almost periodic, it is distal. The converse is not
true; even if $(\Om,\sigma,\R)$ is minimal and distal, it does not need to
be almost periodic. For the main properties of almost periodic and distal
flows we refer the reader to Ellis~\cite{elli}, Sacker and
Sell~\cite{sase}, and Sell~\cite{sell}, \cite{sell2}. The reference
Fink~\cite{fink} describes the basic theory of almost periodic functions
and almost periodic ordinary differential equations.\par
A {\em flow homomorphism\/} from another continuous flow $(Y,\Psi,\R)$ to
$(\Om,\sigma,\R)$ is a continuous map $\pi\colon Y\to \Om$ such that
$\pi(\Psi(t,y))=\sigma(t,\pi(y))$ for every $y\in Y$ and $t\in\R$. If $\pi$
is also bijective, it is called a {\em flow isomorphism\/}. Let $\pi:Y \to
\Om$ be a surjective flow homomorphism and suppose $(Y,\Psi,\R)$ is minimal
(then, so is $(\Om,\sigma,\R)$). $(Y,\Psi,\R)$ is said to be an {\em almost
  automorphic extension\/} of $(\Om,\sigma,\R)$ if there is $\w\in \Om$
such that $\car (\pi^{-1}(\w))=1$. Then, actually $\car (\pi^{-1}(\w))=1$
for $\w$ in a residual subset $\Om_0\subseteq \Om$; in the nontrivial case
$\Om_0\subsetneq \Om$ the dynamics can be very complicated. A minimal flow
$(Y,\Psi,\R)$ is {\em almost automorphic\/} if it is an almost automorphic
extension of an almost periodic minimal flow
$(\Om,\sigma,\R)$. Johnson~\cite{john}, \cite{john2} contain examples of
almost periodic differential equations with almost automorphic solutions
which are not almost periodic. We refer the reader to the work
in~\cite{shyi} for a survey of almost periodic and almost automorphic
dynamics.\par
Let $E$ be a complete metric space and $\R^+=\{t\in\R:t\geq 0\}$. A {\em
  semiflow} $(E,\Phi,\R^+)$ is determined by a continuous map $\Phi:
\R^+\times E \to E,\; (t,x)\mapsto \Phi(t,x)$ which satisfies
\begin{enumerate}
  \renewcommand{\labelenumi}{(\roman{enumi})}
\item $\Phi_0=\text{Id},$
\item $\Phi_{t+s}=\Phi_t \circ \Phi_s\;$ for all $\; t$, $s\in\R^+,$
\end{enumerate}
where $\Phi_t(x)=\Phi(t,x)$ for each $x \in E$ and $t\in \R^+$.  The set
$\{ \Phi_t(x): t\geq 0\}$ is the {\em semiorbit\/} of the point $x$. A
subset $E_1$ of $E$ is {\em positively invariant\/} (or just $\Phi$-{\em
  invariant\/}) if $\Phi_t(E_1)\subset E_1$ for all $t\geq 0$. A semiflow
$(E,\Phi,\R^+)$ admits a {\em flow extension\/} if there exists a
continuous flow $(E,\wt \Phi,\R)$ such that $\wt \Phi(t,x)=\Phi(t,x)$ for
all $x\in E$ and $t\in\R^+$. A compact and positively invariant subset
admits a flow extension if the semiflow restricted to it admits one.\par
Write $\R^-=\{t\in\R\,|\,t\leq 0\}$. A {\em backward orbit\/} of a point
$x\in E$ in the semiflow $(E,\Phi,\R^+)$ is a continuous map $\psi:\R^-\to
E$ such that $\psi(0)=x$ and for each $s\leq 0$ it holds that
$\Phi(t,\psi(s))=\psi(s+t)$ whenever $0\leq t\leq -s$. If for $x\in E$ the
semiorbit $\{\Phi(t,x): t\ge 0\}$ is relatively compact, we can consider
the {\em omega-limit set\/} of $x$,
\[
\mathcal{O}(x)=\bigcap_{s\ge 0}{\rm cls}{\{\Phi(t+s,x): t\ge 0\}}\,,
\]
which is a nonempty compact connected and $\Phi$-invariant set.  Namely, it
consists of the points $y\in E$ such that $y=\lim_{\,n\to \infty}
\Phi(t_n,x)$ for some sequence $t_n\uparrow \infty$. It is well-known that
every $y\in\mathcal{O}(x)$ admits a backward orbit inside this
set. Actually, a compact positively invariant set $M$ admits a flow
extension if every point in $M$ admits a unique backward orbit which
remains inside the set $M$ (see~\cite{shyi}, part~II).\par
A compact positively invariant set $M$ for the semiflow $(E,\Phi,\R^+)$ is
{\em minimal\/} if it does not contain any other nonempty compact
positively invariant set than itself. If $E$ is minimal, we say that the
semiflow is minimal.\par
A semiflow is {\em of skew-product type\/} when it is defined on a vector
bundle and has a triangular structure; more precisely, a semiflow
$(\Om\times X,\tau,\,\R^+)$ is a {\em skew-product\/} semiflow over the
product space $\Om\times X$, for a compact metric space $(\Om,d)$ and a
complete metric space $(X,\textsf{d})$, if the continuous map $\tau$ is as
follows:
\begin{equation}\label{skewp}
  \begin{array}{cccl}
    \tau \colon  &\R^+\times\Om\times X& \longrightarrow & \Om\times X \\
    & (t,\w,x) & \mapsto &(\w{\cdot}t,u(t,\w,x))\,,
  \end{array}
\end{equation}
where $(\Om,\sigma,\R)$ is a real continuous flow
$\sigma:\R\times\Om\rightarrow\Om$, $\,(t,\w)\mapsto \w{\cdot}t$, called
the {\em base flow\/}. The skew-product semiflow~\eqref{skewp} is {\em
  linear\/} if $u(t,\w,x)$ is linear in $x$ for each
$(t,\w)\in\R^+\times\Om$. The definitions of stability, uniform stability
and asymptotic stability for non-autonomous differential equations used
along the paper can be found in~\cite{sell}, \cite{sell2} (see also Conley
and Miller~\cite{comi} for an interesting remark).
\section{Non-Autonomous stable linear
  $D$-operators}\label{section_operator}
Given an $m\times m$ matrix $\mu=[\mu_{ij}]_{ij}$ of measures with finite
total variation on a measurable space $(Y,\zeta)$ and a measurable subset
of $Y$, $E\in\zeta$, $|\mu_{ij}|(E)$ will denote the total variation of
$\mu_{ij}$ over $E$; besides, the maximum norm of the $m\times m$ matrix
$[|\mu_{ij}|(E)]_{ij}$ will be denoted by $\n{\mu}_\infty(E)$ and the
$m\times m$ matrix of positive measures $[|\mu_{ij}|]_{ij}$ will be denoted
by $|\mu|$.\par
Let $(\Om,d)$ be a compact metric space and let $\sigma:\R\times\Om\to\Om$
be a continuous real flow on $\Om$. We will denote
$\w{\cdot}t=\sigma(\w,t)$, $t\in\R$, $\w\in\Om$. We will assume in the
remainder of the paper that the flow $\sigma$ is minimal.\par
Let $X=C((-\infty,0],\R^m)$, which is a Fr{\'e}chet space when endowed with
the compact-open topology, i.e.~the topology of uniform convergence over
compact subsets. This topology happens to be metric for the distance
\[
\di(x,y)=\sum_{n=1}^\infty \frac{1}{2^n}\frac
{\n{x-y}_n}{1+\n{x-y}_n}\,,\quad x,y\in X\,,
\]
where $\n{x}_n=\sup_{s\in[-n,0]}\n{{x(s)}}$, and $\n{\cdot}$ denotes the
maximum norm in $\R^m$. Let $BU\subset X$ be the Banach space
\[
BU=\{x\in X:x \text{ is bounded and uniformly continuous}\}
\]
with the supremum norm $\n{x}_\infty=\sup_{s\in(-\infty,0]} \n{x(s)}$.
Given $r>0$, we will denote
\[
B_r=\{x\in BU:\n{x}_\infty \leq r\}.
\]
As usual, given $I=(-\infty,a]\subset\R$, $t\in I$ and a continuous
function $x:I\to\R^m$, $x_t$ will denote the element of $X$ defined by
$x_t(s)=x(t+s)$ for $s\in (-\infty,0]$.\par
Let $D:\Om\times BU\to\R^m$ be an operator satisfying the following
hypotheses:
\begin{itemize}
\item[(D1)] $D$ is linear and continuous in its second variable and the map
  $\Om\to\mathcal{L}(BU,\R^m)$, $\w\mapsto D(\w,{\cdot})$ is continuous;
\item[(D2)] for each $r>0$, $D\colon\Om\times B_r\to \R^m$ is continuous
  when we take the restriction of the compact-open topology to $B_r$,
  i.e. if $\w_n\to\w$ and $x_n\stackrel{\di}\to x$ as $n\to\infty$ with
  $(\w_n,x_n)$, $(\w,x)\in\Om\times B_r$, then
  $\lim_{n\to\infty}D(\w_n,x_n)=D(\w,x)$.
\end{itemize}
\begin{lemma}\label{D_measure}
  For each $\w\in\Om$, there exists an $m\times m$ matrix
  $\mu(\w)=[\mu_{ij}(\w)]_{ij}$ of real Borel regular measures with finite
  total variation such that
  \[
  D(\w,x)=\int_{-\infty}^0[d\mu(\w)]x,\quad(\w,x)\in\Om\times BU.
  \]
\end{lemma}
\begin{proof}
  This result follows by applying Proposition 3.1 in~\cite{muno} to each
  $\w\in\Om$.
\end{proof}
For each $\w\in\Om$, let $B(\w)=\mu(\w)(\{0\})$. Now let
$\nu(\w)=B(\w)\delta_0-\mu(\w)$, $\w\in\Om$, where $\delta_0$ is the Dirac
measure at 0, that is, $\int_{-\infty}^0[d\delta_0]x=x(0)$ for all $x\in
BU$. It is clear that $|\nu_{ij}(\w)|(\{0\})=0$ for all
$i,\,j\in\{1,\ldots,m\}$ and all $\w\in\Om$. Besides, from the dominated
convergence theorem, it follows that
\[
\lim_{\rho\to 0^+}|\nu_{ij}(\w)|([-\rho,0])=0\text{ and }\lim_{\rho\to
  \infty}|\nu_{ij}(\w)|((-\infty,-\rho])=0
\]
for each $\w\in\Om$.
\begin{proposition}\label{B_cont}
  $B:\Om\to\mathbb{M}_m(\R)$, $\w\mapsto B(\w)$ is a continuous map.
\end{proposition}
\begin{proof}
  For each $\rho>0$, let $\varphi_\rho:(-\infty,0]\to\R$ be the function
  given for $s\leq 0$ by
  \begin{equation}\label{phi_r}
    \varphi_\rho(s)=\left\{
      \begin{array}{ll}
        0&\text{if }s\leq -2\rho,\\
        \rho^{-1}s+2&\text{if }-2\rho<s\leq -\rho,\\
        1&\text{if }-\rho<s\leq 0.
      \end{array}
    \right.
  \end{equation}
  Let $i\in\{1,\ldots,m\}$ and let $\{\w_n\}_n\subset\Om$ be a sequence
  converging to some $\w_0\in\Om$. On the one hand, a straightforward
  application of the dominated convergence theorem yields
  \[
  \lim_{\rho\to 0^+}\int_{-\infty}^0[d\mu(\w_n)]\varphi_\rho
  e_i=B(\w_n)e_i\text{ for all }n\in\N\text{ and }\lim_{\rho\to
    0^+}D(\w_0,\varphi_\rho e_i)=B(\w_0)e_i.
  \]
  On the other hand, from (D1) we deduce that
  \[
  \lim_{n\to\infty}\int_{-\infty}^0[d\mu(\w_n)]\varphi_\rho
  e_i=D(\w_0,\varphi_\rho e_i)
  \]
  uniformly for $\rho>0$, and the result follows immediately.
\end{proof}
\begin{proposition}\label{L_cont}
  Let $L:\Om\times BU\to\R^m$, $(\w,x)\mapsto B(\w)x(0)-D(\w,x)$. Then the
  mapping $\Om\to\mathcal{L}(BU,\R^m)$, $\w\mapsto L(\w,{\cdot})$ is
  continuous. Equivalently, for every sequence $\{\w_n\}_n\subset\Om$
  converging to $\w_0\in\Om$ and all $i,\,j\in\{1,\ldots,m\}$, we have that
  $\lim_{n\to\infty}|\nu_{ij}(\w_n)-\nu_{ij}(\w_0)|((-\infty,0])=0$.
\end{proposition}
\begin{proof}
  Let $\w_1,\,\w_2\in\Om$ and $x\in BU$ with $\n{x}_\infty\leq 1$; then
  \[
  \begin{split}
    \n{L(\w_1,x&)-L(\w_2,x)}\leq\n{B(\w_1)-B(\w_2)}+\n{D(\w_1,x)-D(\w_2,x)}.
  \end{split}
  \]
  The result follows from (D1) and Proposition~\ref{B_cont}.
\end{proof}
\begin{corollary}\label{nu_unif_to_0}
  Under hypotheses {\upshape (D1)--(D2)}, the following statements hold:
  \begin{itemize}
  \item[(i)] $\lim_{\rho\to 0^+}\n{\nu(\w)}_\infty([-\rho,0])=0$ uniformly
    for $\w\in\Om$;
  \item[(ii)] $\lim_{\rho\to\infty}\n{\nu(\w)}_\infty((-\infty,-\rho])=0$
    uniformly for $\w\in\Om$.
  \end{itemize}
\end{corollary}
\begin{proof}
  Let us prove (i). Note that, for each $\w_1,\,\w_2\in\Om$ and each
  $i,\,j\in\{1,\ldots,m\}$, we have that
  \[
  \left||\nu_{ij}(\w_1)|([-\rho,0])-|\nu_{ij}(\w_2)|([-\rho,0])\right|\leq
  |\nu_{ij}(\w_1)-\nu_{ij}(\w_2)|((-\infty,0]).
  \]
  Moreover, $\lim_{\rho\to 0^+}\n{\nu(\w)}_\infty([-\rho,0])=0$ for all
  $\w\in\Om$. This is a family of continuous functions decreasing to 0. As
  a result, Dini's theorem implies that this family converges to 0
  uniformly for $\w\in\Om$. The proof of (ii) is analogous.
\end{proof}
Let us assume one more hypothesis on the operator $D$, which is a natural
generalization of the atomic character as seen in~\cite{hale}
and~\cite{have}:
\begin{itemize}
\item[(D3)] $B(\w)$ is a regular matrix for all $\w\in\Om$.
\end{itemize}
\begin{theorem}\label{existence}
  For all $h\in C([0,\infty),\R^m)$ and all $(\w,\varphi)\in\Om\times BU$
  with $D(\w,\varphi)=h(0)$, there exists $x\in C(\R,\R^m)$ such that
  \addtocounter{equation}{1}
  \begin{equation}\tag*{(\arabic{section}.\arabic{equation})$_\w$}\label{nonhomogeneous}
    \left\{
      \begin{array}{ll}
        D(\w{\cdot}t,x_t)=h(t),&t\geq 0,\\
        x_0=\varphi.
      \end{array}
    \right.
  \end{equation}
  Moreover, given $\rho>0$, there are positive constants $k_\rho^1$,
  $k_\rho^2$ such that for each $t\in[0,\rho]$
  \[
  \n{x_t}_\infty\leq k_\rho^1\,\sup_{0\leq u\leq t}\n{h(u)}+
  k_\rho^2\,\n{\varphi}_\infty\,.
  \]
\end{theorem}
This bound for the solution leads us to its uniqueness. Namely, if
$x^1,\,x^2$ are solutions of the equation, then for all $t\geq 0$, we have
\[
\left\{
  \begin{array}{ll}
    D(\w{\cdot}t,x^1_t-x^2_t)=0\,, & t\geq 0\,, \\
    x^1_0-x^2_0=0,
  \end{array} \right.
\]
Thus, given $t>0$, $\n{x^1(t)-x^2(t)}\leq k^1_t\,0+k^2_t\,0=0$.
\begin{theorem}\label{Dhat_cont}
  Let us define the map
  \begin{equation}\label{Dhat}
    \begin{array}{lcclrcl}
      \widehat{D}:&\Om\times BU&\longrightarrow & \Om\times BU &&\\
      & (\w,x) & \mapsto &(\w,\widehat{D}_2(\w,x))
    \end{array}
  \end{equation}
  where $\widehat D_2(\w,x):(-\infty,0]\to\R^m$, $s\mapsto
  D(\w{\cdot}s,x_s)$. Then $\widehat D$ is well defined, $\widehat D_2$ is
  linear and continuous for the norm in its second variable for all
  $\w\in\Om$ and, for all $r>0$, $\widehat D$ is uniformly continuous on
  $\Om\times B_r$ when we take the restriction of the compact-open topology
  to $B_r$. In addition, if the flow $(\Om,\sigma,\R)$ is almost periodic,
  then $\widehat D$ is uniformly continuous on $\Om\times B_r$ for all
  $r>0$ when we take the norm on $B_r$.
\end{theorem}
\begin{proof}
  Let us check that $\widehat D$ is well defined. Let $(\w,x)\in\Om\times
  BU$ and let $h=\widehat{D}_2(\w,x):(-\infty,0]\to\R^m$. From (D1) and the
  uniform continuity of $\sigma$ on, say, $[0,1]\times\Om$, it follows
  that, for all $\eps>0$, there exists $\delta\in(0,1)$ such that, if
  $t,\,s\leq 0$ and $|t-s|<\delta$ then
  \[
  \n{D(\w{\cdot}t,{\cdot})-D(\w{\cdot}s,{\cdot})}\,\n{x}_\infty\leq\frac{\eps}{2}\,\,\text{
    and
  }\,\sup_{\w_1\in\Om}\n{D(\w_1,{\cdot})}\,\n{x_t-x_s}_\infty\leq\frac{\eps}{2},
  \]
  whence
  \[
  \n{D(\w{\cdot}t,x_t)-D(\w{\cdot}s,x_s)}\leq\eps.
  \]
  Clearly,
  $\n{h}_\infty\leq\sup_{\w_1\in\Om}\n{D(\w_1,{\cdot})}\,\n{x}_\infty$ and,
  consequently, $h\in BU$. This way, $\wh D$ is well defined.\par
  The linearity of $\widehat D_2$ in its second variable is clear. Besides,
  for all $\w\in\Om$, the continuity of $\widehat D_2(\w,{\cdot})$ for the
  norm on $BU$ is a straightforward consequence of (D1).\par
  Let us check the uniform continuity of $\widehat D$ on $\Om\times B_r$,
  $r\geq 0$, when we take the restriction of the compact-open topology to
  $B_r$. In order to do so, let us fix $\rho>0$ and $\eps>0$; there is a
  $\delta>0$ such that, for all $\w^1,\,\w^2\in\Om$ with
  $d(\w^1,\w^2)<\delta$ and all $s\in[-\rho,0]$, it holds that
  $\n{D(\w^1{\cdot}s,{\cdot})-D(\w^2{\cdot}s,{\cdot})}<\eps/(2r)$. Thanks
  to Corollary~\ref{nu_unif_to_0}, there exists $\rho_0>0$ such that
  $\sup_{\w\in\Om}\n{\mu(\w)}_\infty((-\infty,-\rho_0])<\eps/(8r)$. Now,
  let $(\w^1,x^1)$, $(\w^2,x^2)\in\Om\times B_r$ such that
  $d(\w^1,\w^2)<\delta$ and satisfying
  $\sup_{\w\in\Om}\n{\mu(\w)}_\infty((-\infty,0])
  \n{x^1-x^2}_{[-\rho-\rho_0,0]}<\eps/4$. If $s\in[-\rho,0]$, then
  \[
  \begin{split}
    \n{D(\w^1{\cdot}s,x^1_s)&-D(\w^2{\cdot}s,x^2_s)}\leq\n{D(\w^1{\cdot}s,{\cdot})-D(\w^2{\cdot}s,{\cdot})}\,r+\n{D(\w^2{\cdot}s,(x^1-x^2)_s)}\\
    \leq&\frac{\eps}{2r}r+\left\|\int_{-\infty}^{-\rho_0}[d\mu(\w^2{\cdot}s)](x^1-x^2)_s\right\|
    +\left\|\int_{-\rho_0}^0[d\mu(\w^2{\cdot}s)](x^1-x^2)_s\right\|\\
    \leq&\frac{\eps}{2}+\frac{\eps}{8r}2r+\frac{\eps}{4}=\eps.
  \end{split}
  \]
  This inequality yields the expected result.\par
  Finally, assume that $(\Om,\sigma,\R)$ is almost periodic. Let us prove
  that $\widehat D_2$ is uniformly continuous on $\Om\times B_r$, $r>0$,
  when we take the norm on $B_r$. From (D1) and the almost periodicity of
  $(\Om,\sigma,\R)$, it follows that, for all $\eps>0$, there exists
  $\delta>0$ such that, if $\w^1,\,\w^2\in\Om$, $s\leq 0$ and
  $d(\w^1,\w^2)<\delta$, then
  $\n{D(\w^1{\cdot}s,{\cdot})-D(\w^2{\cdot}s,{\cdot})}<\eps/(2r)$. Taking
  $(\w^1,x^1)$, $(\w^2,x^2)\in\Om\times B_r$ with $d(\w^1,\w^2)<\delta$ and
  such that $\sup_{\w\in\Om}\n{\mu(\w)}_\infty((-\infty,0])
  \n{x^1-x^2}_\infty<\eps/2$, an argument similar to the previous one
  yields the desired property.
\end{proof}
\begin{definition}
  The mapping $D$ is said to be {\em stable\/} if there is a continuous
  function $c\in C([0,\infty),\R)$ with $\lim_{\,t\to\infty}c(t)=0$ such
  that, for each $(\w,\varphi)\in\Om\times BU$ with $D(\w,\varphi)=0$, the
  solution of the homogeneous problem
  \[\left\{
    \begin{array}{ll}
      D(\w{\cdot}t,x_t)=0\,, & t\geq 0 \\
      x_0=\varphi\,,
    \end{array} \right.
  \] satisfies $\n{x(t)}\leq c(t)\,\n{\varphi}_\infty$ for each $t\geq 0$.
\end{definition}
The following statement, whose proof is similar to the one in~\cite{hale},
provides a nonhomogeneous version of the concept of stability for a
$D$-operator.
\begin{theorem}\label{D_stable}
  Let us assume that $D$ is stable. Then there are a continuous function
  $c\in C([0,\infty),\R)$ with $\lim_{\,t\to\infty}c(t)=0$ and a positive
  constant $k>0$ such that the solution of the equation
  \[
  \left\{
    \begin{array}{ll}
      D(\w{\cdot}t,x_t)=h(t)\,, & t\geq 0\,, \\
      x_0=\varphi\,,
    \end{array} \right.
  \]
  where $h\in C([0,\infty),\R^m)$, $(\w,\varphi)\in\Om\times BU$ and
  $D(\w,\varphi)=h(0)$, satisfies
  \[\n{x(t)}\leq c(t)\,\n{\varphi}_\infty+ k\,\sup_{0\leq
    u\leq t}\n{h(u)}\] for each $t\geq 0$.
\end{theorem}
It is easy to see that $\widehat D_2$ has the following integral
representation:
\[
\widehat{D}_2(\w,x)(s)=B(\w{\cdot}s)x(s)-\int_{-\infty}^0[d\nu(\w{\cdot}s)]\,x_s
\]
for each $(\w,x)\in\Om\times BU$ and $s\in(-\infty,0]$.\par
We are now in a position to state the main theorem of this section which
assures the invertibility of $\wh D$ and specifies the regularity of its
inverse when the linear operator $D$ is stable.
\begin{theorem}\label{D-1}
  Under hypotheses {\upshape(D1)--(D3)}, the following statements hold:
  \begin{itemize}
  \item[(i)] let us assume that $D$ is stable; then there is a positive
    constant $k>0$ such that $\n{x^h}_\infty\leq k\,\n{h}_\infty$ for all
    $h\in BU,\,\w\in\Om$ and $x^h\in BU$ satisfying
    $D(\w{\cdot}s,x_s^h)=h(s)$ for $s\leq 0$;
  \item[(ii)] if $D$ is stable, then $\widehat{D}$ is invertible,
    $\widehat{D}^{-1}_2(\w,{\cdot})$ is linear and continuous on $BU$ for
    all $\w\in\Om$ when we consider the norm on $BU$, and
    $\widehat{D}^{-1}$ is uniformly continuous on $\Om\times B_r$ for all
    $r>0$ when we take the restriction of the compact-open topology to
    $B_r$. In addition, if the flow $(\Om,\sigma,\R)$ is almost periodic,
    $\widehat{D}^{-1}$ is uniformly continuous on $\Om\times B_r$ for all
    $r>0$ when we take the norm on $B_r$;
  \item[(iii)] let us assume that, for each $r>0$ and each sequence
    $\{(\w_n,x_n)\}_n\subset\Om\times BU$ such that $\n{\widehat D_2
      (\w_n,x_n)}_\infty\leq r$, $\w_n\to\w\in\Om$ and {\upshape$\widehat
      D_2(\w_n,x_n) \stackrel{\di}\to 0$} as $n\to\infty$, it holds that
    $x_n(0)\to 0$ as $n\to\infty$; then $D$ is stable.
  \end{itemize}
\end{theorem}
\begin{proof}
  The proof of (i) is analogous to the one of Proposition~3.7
  in~\cite{muno} due to the uniform character over $\Om$ of the function
  $c$ and the constant $k$ of Theorem~\ref{D_stable}.\par
  Now, we will prove statement (ii). $\widehat{D}$ is injective because, if
  we have $(\w^1,x^1)$, $(\w^2,x^2)\in\Om\times BU$ with $\widehat
  D(\w^1,x^1)=\widehat D(\w^2,x^2)$, then $\w^1=\w^2$ and, from (i) and the
  fact that $D(\w^1{\cdot}s,x^1_s-x^2_s)=0$ for $s\leq 0$, we get
  $x^1=x^2$.\par
  In order to show that $\widehat{D}$ is surjective, let
  $(\w,h)\in\Om\times BU$ and $\{h_n\}_n\subseteq B_r$, for some $r>0$, be
  a sequence of continuous functions whose components are of compact
  support such that $h_n\stackrel{\textsf{d}\;}\to h$ as
  $n\to\infty$. Moreover, it is easy to choose them with the same modulus
  of uniform continuity as $h$. Let us check that, for each $n\in\N$, there
  is an $x^n\in BU$ such that $\widehat{D}_2(\w,x^n)=h_n$, that is,
  $D(\w{\cdot}s,x_s^n)=h_n(s)$ for $s\leq 0$ and $n\in\N$. Fix $n\in\N$ and
  $\rho_n>0$ such that $\text{supp}(h_n)\subset[-\rho_n,0]$. Let $\wt
  h_n:[0,\infty)\to\R^m$ be the function defined for $t\geq 0$ by
  \[
  \wt h_n(t)=
  \begin{cases}
    h_n(t-\rho_n)&\text{if }t\in[0,\rho_n],\\
    h_n(0)&\text{if }t\geq \rho_n.
  \end{cases}
  \]
  Since $\wt h_n(0)=0$, by Theorem~\ref{existence}, we have that there
  exists $\wt x^n\in C(\R,\R^m)$ such that $D(\w{\cdot}(t-\rho_n),\wt
  x_t^n)=\wt h_n(t)$ for $t\geq 0$ and $\wt x_0^n=0$. Let
  $x^n:(-\infty,0]\to\R^m$ be the function defined by $x^n(s)=\wt
  x^n(s+\rho_n)$ for $s\leq 0$. Clearly, the function $x^n$ is continuous
  and of compact support. Now, if $s\in[-\rho_n,0]$, then
  $D(\w{\cdot}s,x_s^n)=D(\w{\cdot}(-\rho_n+(s+\rho_n)),\wt
  x_{s+\rho_n}^n)=\wt h_n(s+\rho_n)=h_n(s)$ and, if $s\leq -\rho_n$, then
  $D(\w{\cdot}s,x_s^n)=D(\w{\cdot}(-\rho_n+(s+\rho_n)),\wt
  x_{s+\rho_n}^n)=D(\w{\cdot}s,0)=h_n(s)=0$ as wanted.\par
  From (i), there exists $k>0$ such that $\n{x^n}_\infty\leq
  k\n{h_n}_\infty\leq kr$. Let us fix $\eps>0$; since the restriction of
  $\sigma$ to $[-1,0]\times\Om$ is uniformly continuous, we can fix
  $\delta>0$ such that, for all $\tau\in[-\delta,0]$ and all $s\leq 0$,
  $\n{h_n-(h_n)_\tau}_\infty<\eps/(2k)$ and
  $\n{D(\w{\cdot}(s+\tau),{\cdot})-D(\w{\cdot}s,{\cdot})}<\eps/(2k^2r)$. For
  each $n\in\N$ and each $\tau\in[-\delta,0]$, let
  $g^\tau_n:(-\infty,0]\to\R^m$, $s\mapsto
  D(\w{\cdot}s,(x^n-x^n_\tau)_s)$. Then, for all $s\leq 0$, all
  $\tau\in[-\delta,0]$ and all $n\in\N$,
  \[
  \begin{split}
    \n{g^\tau_n(s)}\leq&\n{D(\w{\cdot}s,x^n_s)\!-\!D(\w{\cdot}(s+\tau),x^n_{s+\tau})}+\n{D(\w{\cdot}(s+\tau),x^n_{s+\tau})\!-\!D(\w{\cdot}s,x^n_{s+\tau})}\\
    \leq&\frac{\eps}{2k}+\frac{\eps}{2k^2r}\n{x^n}_\infty\leq\frac{\eps}{k}.
  \end{split}
  \]
  From (i), we deduce again that
  \[
  \n{x^n-x^n_\tau}_\infty\leq k\,\n{g^\tau_n}_\infty\leq
  k\,\frac{\eps}{k}=\eps
  \]
  for all $n\in\N$ and all $\tau\in[-\delta,0]$. Thus $\{x^n\}_n$ is
  equicontinuous and, consequently, relatively compact for the compact-open
  topology. Hence, there is a convergent subsequence of $\{x^n\}_n$, let us
  assume the whole sequence, i.e. there is a continuous function $x$ such
  that $x^n\stackrel{\di}\to x$ as $n\to\infty$. Therefore, we have that
  $\n{x}_\infty\leq k\,r$ and $\n{x^n(s)-x^n(s+t)}\to\n{x(s)-x(s+t)}$ as $n\to\infty$
  for all $s,t\leq 0$, which implies that $x\in BU$. From this, $x^n_s
  \stackrel{\textsf{d}\;}\to x_s$ for each $s\leq 0$ and the expression of
  $D$ yields $D(\w{\cdot}s,x_s^n)=h_n(s)\to D(\w{\cdot}s,x_s)$, i.e.
  $D(\w{\cdot}s,x_s)=h(s)$ for $s\leq 0$ and $\widehat D_2(\w,x)=h$. Then
  $\widehat{D}$ is surjective, as claimed.\par
  Let us check that $\widehat{D}^{-1}$ is uniformly continuous on
  $\Om\times B_r$, $r>0$, when we take the restriction of the compact-open
  topology to $B_r$. Fix $\eps>0$ and $\rho>0$; using
  Theorem~\ref{D_stable}, it is clear that we can find a $\rho_0>0$ such
  that $c(t)<\eps/(4kr)$ for all $t\geq \rho_0$. Besides, there is a
  $\delta>0$ such that, if $\w_1,\,\w_2\in\Om$ and $d(\w_1,\w_2)<\delta$,
  then $\n{D(\w_1{\cdot}s,x)- D(\w_2{\cdot}s,x)}<\eps/(4k)$ for all
  $s\in[-\rho_0-\rho,0]$ and all $x\in B_{kr}$, thanks to the uniform
  continuity of $\sigma$ on $[-\rho_0-\rho,0]\times\Om$ and (D1). Let
  $(\w_1,h_1),\,(\w_2,h_2)\in\Om\times B_r$ such that $d(\w_1,\w_2)<\delta$
  and $\n{h_1-h_2}_{[-\rho_0-\rho,0]}\leq\eps/(4k)$;
  $x_1=\widehat{D}^{-1}_2(\w_1,h_1)$, $x_2=\widehat{D}^{-1}_2(\w_2,h_2)\in
  B_{kr}$ thanks to (i). Then, for all $t\in[0,\rho_0+\rho]$,
  \[
  D(\w_i{\cdot}(-\rho_0-\rho+t),(x_i)_{-\rho_0-\rho+t})=h_i(-\rho_0-\rho+t),\quad
  i=1,2,
  \]
  whence it follows that, for all $t\in[0,\rho_0+\rho]$,
  \[
  \begin{split}
    \n{D(\w_1{\cdot}&(-\rho_0-\rho+t),(x_1)_{-\rho_0-\rho+t})-D(\w_1{\cdot}(-\rho_0-\rho+t),(x_2)_{-\rho_0-\rho+t})}\leq\\
    \leq&\n{D(\w_1{\cdot}(-\rho_0-\rho+t),(x_1)_{-\rho_0-\rho+t})-D(\w_2{\cdot}(-\rho_0-\rho+t),(x_2)_{-\rho_0-\rho+t})}\\
    &+\n{D(\w_2{\cdot}(-\rho_0-\rho+t),(x_2)_{-\rho_0-\rho+t})-D(\w_1{\cdot}(-\rho_0-\rho+t),(x_2)_{-\rho_0-\rho+t})}\\
    \leq&\n{h_1(-\rho_0-\rho+t)-h_2(-\rho_0-\rho+t)}+\frac{\eps}{4k}\leq\frac{\eps}{4k}+\frac{\eps}{4k}=\frac{\eps}{2k}.
  \end{split}
  \]
  Now, the stability of $D$ and Theorem~\ref{D_stable} yield that, for all
  $t\in[\rho_0,\rho_0+\rho]$,
  \[
  \n{x_1(-\rho_0-\rho+t)-x_2(-\rho_0-\rho+t)}\leq
  \!\!\sup_{t\in[\rho_0,\rho_0+\rho]}\!\!c(t)\n{x_1-x_2}_\infty+k\frac{\eps}{2k}\leq\frac{\eps}{4kr}2kr+\frac{\eps}{2}=\eps,
  \]
  that is, $\n{x_1-x_2}_{[-\rho,0]}\leq\eps$ and the expected result
  holds.\par
  Finally, we prove that, provided that $(\Om,\sigma,\R)$ is almost
  periodic, $\widehat{D}^{-1}$ is uniformly continuous on $\Om\times B_r$,
  $r>0$, when the norm is considered on $B_r$. Let $\eps>0$; from
  Theorem~\ref{Dhat_cont} and the almost periodicity of $(\Om,\sigma,\R)$,
  it follows that there exists $\delta>0$ such that, for all
  $\w_1,\,\w_2\in\Om$ with $d(\w_1,\w_2)<\delta$ and all $x\in B_{kr}$,
  $\n{\wh D_2(\w_1,x)-\wh D_2(\w_2,x)}_\infty<\eps/(2k)$. Now, fix
  $(\w_1,h_1),\,(\w_2,h_2)\in\Om\times B_r$ with $d(\w_1,\w_2)<\delta$ and
  $\n{h_1-h_2}_\infty\leq\eps/(2k)$. From (i), we have that $x_1=(\wh
  D^{-1})_2(\w_1,h_1),\,x_2=(\wh D^{-1})_2(\w_2,h_2)\in B_{kr}$. Let
  $y=(\wh D^{-1})_2(\w_1,h_2)$; using (i),
  \[
  \frac{1}{k}\n{y-x_2}_\infty\leq\n{\wh D_2(\w_1,y)-\wh
    D_2(\w_1,x_2)}_\infty=\n{\wh D_2(\w_2,x_2)-\wh
    D_2(\w_1,x_2)}_\infty\leq\frac{\eps}{2k},
  \]
  which, together with (i) again, yields
  \[
  \n{x_1-x_2}_\infty\leq\n{x_1-(\wh D^{-1})_2(\w_1,h_2)}_\infty+\n{(\wh
    D^{-1})_2(\w_1,h_2)-x_2}_\infty\leq
  k\frac{\eps}{2k}+\frac{\eps}{2}=\eps,
  \]
  as desired.\par
  Now we prove statement (iii). Let $C_D=\{(\w,\varphi)\in\Om\times BU:
  D(\w,\varphi)=0\}$. For each $\rho>0$, we define
  $\mathcal{L}_\rho:C_D\to\R^m$, $(\w,\varphi)\mapsto x(\rho)$, where $x$
  is the solution of
  \[\left\{
    \begin{array}{ll}
      D(\w{\cdot}t,x_t)=0\,, & t\geq 0\,, \\
      x_0=\varphi\,.
    \end{array} \right.
  \]
  Let us observe that $C_D=\{(\w,\varphi):\w\in\Om,\,\varphi\in C_D(\w)\}$,
  where $C_D(\w)=\{\varphi\in BU:D(\w,\varphi)=0\}$ is a vector space for
  each $\w\in\Om$.\par
  From the uniqueness of the solution of~\ref{nonhomogeneous}, it is easy
  to check that $\mathcal{L}_\rho$ is well defined and linear in its second
  variable. In addition, from Theorem~\ref{existence}, we deduce that
  $\n{\mathcal{L}_\rho(\w,\varphi)}=\n{x(\rho)}\leq
  k_\rho^2\,\n{\varphi}_\infty$ for all $(\w,\varphi)\in C_D$, whence
  $\n{\mathcal L_\rho(\w,{\cdot})}\leq k_\rho^2$ for all $\w\in\Om$.\par
  Next, we check that
  $\sup_{\w\in\Om}\n{\mathcal{L}_\rho(\w,{\cdot})}_\infty\to 0$ as
  $\rho\to\infty$, which shows the stability of $D$ because
  $\n{x(\rho)}\leq c(\rho)\,\n{\varphi}_\infty$ for all $(\w,\varphi)\in
  C_D$, where
  $c(\rho)=\sup_{\w\in\Om}\n{\mathcal{L}_\rho(\w,{\cdot})}_\infty$. Let us
  assume, on the contrary, that there exist $\delta>0$, a sequence
  $\rho_n\uparrow\infty$ and a sequence $\{\varphi_n\}_n$ with
  $\varphi_n\in C_D(\w_n)$, $\n{\varphi_n}_\infty\leq 1$ and such that
  $\n{\mathcal{L}_{\rho_n}(\w_n,\varphi_n)}\geq \delta$ for each $n\in\N$.
  That is, $\n{x^n(\rho_n)}\geq \delta$ where $x^n$ is the solution of
  \[\left\{
    \begin{array}{ll}
      D(\w_n{\cdot}t,x^n_t)=0\,, & t\geq 0\,, \\
      x_0^n=\varphi_n\,.
    \end{array} \right.
  \]
  Therefore,
  \[
  \begin{cases}
    D(\w_n{\cdot}(\rho_n+s),(x_{\rho_n}^n)_s)=D(\w_n{\cdot}(\rho_n+s),x_{\rho_n+s}^n)=0 & \text{ if } s\in[-\rho_n,0]\,, \\
    D(\w_n{\cdot}(\rho_n+s),(x_{\rho_n}^n)_s)=D(\w_n{\cdot}(\rho_n+s),(\varphi_n)_{\rho_n+s})
    & \text{ if } s\leq -\rho_n\,,
  \end{cases}
  \]
  and taking $r=\sup_{\w\in\Om}\n{D(\w,{\cdot})}$, the sequence
  $\{x_{\rho_n}^n\}_{n\in\N}\subset BU$ satisfies that $\n{\widehat
    D(\w_n{\cdot}\rho_n, x^n_{\rho_n})}_\infty\leq r$ and $\widehat
  D(\w_n{\cdot}\rho_n, x_{\rho_n}^n)\stackrel{\di}\to 0$ as $n\to
  \infty$. We can assume without loss of generality that
  $\w_n{\cdot}\rho_n\to\w\in\Om$. Consequently,
  $x_{\rho_n}^n(0)=x^n(\rho_n)\to 0$ as $n\to \infty$, which contradicts
  the fact that $\n{x^n(\rho_n)}\geq \delta$, and finishes the proof.
\end{proof}
As a consequence, the operator $D$ is stable if and only if $\wh D$ is
invertible and $\wh D^{-1}$ is uniformly continuous on $\Om\times B_r$ for
the compact-open topology for all $r>0$.\par
The following statement provides a symmetric theory for the operators $\wh
D$ and $\wh D^{-1}$. In particular, $\wh D^{-1}$ is generated by a linear
operator $D^*$ which satisfies (D1)--(D3). We omit the proof, which follows
the arguments of Proposition 2.10 in~\cite{muno}.
\begin{proposition}\label{D-1_like_D}
  Suppose that $D$ is stable and define
  \[
  \begin{array}{lccl}
    D^*: &\Om\times BU &\longrightarrow & \R^m \\
    & (\w,x) & \mapsto &(\widehat{D}^{-1})_2(\w,x)(0).
  \end{array}
  \]
  Then $D^*$ also satisfies {\upshape (D1)--(D3)} and is stable. Moreover,
  for all $s\leq 0$ and all $(\w,x)\in\Om\times BU$, it holds that
  $(\widehat{D}^{-1})_2(\w,x)(s)=D^*(\w{\cdot}s,x_s)$. In particular, for
  all $\w\in\Om$, there is an $m\times m$ matrix
  $\mu^*(\w)=[\mu^*_{ij}(\w)]_{ij}$ of real Borel regular measures with
  finite total variation such that
  \[
  (\widehat{D}^{-1})_2(\w,x)(s)=\int_{-\infty}^0[d\mu^*(\w{\cdot}s)]x_s,
  \quad(\w,x)\in\Om\times BU,\,s\leq 0.
  \]
\end{proposition}
\section{Transformed exponential order and structure of omega-limit
  sets}\label{section_copy_base}
Let $F:\Om\times BU\to \R^m$, let $(\Om,d)$ be a compact metric space and
let $\R\times\Om\to\Om$, $\w\mapsto\w{\cdot}t$ be a minimal real flow on
$\Om$. Let $D:\Om\times BU\to\R^m$ be a stable operator satisfying
hypotheses (D1)--(D3).\par
Let us consider the family of equations \addtocounter{equation}{1}
\begin{equation}\tag*{(\arabic{section}.\arabic{equation})$_\w$}\label{original_family}
  \frac{d}{dt}D(\w{\cdot}t,z_t)=F(\w{\cdot}t,z_t),\quad t\geq 0,\,\w\in\Om.
\end{equation}
We take the componentwise partial order relation on $\R^m$, that is, if
$v,\,w\in\R^m$, then
\[
\begin{split}
  v\le w \quad&\Longleftrightarrow\quad v_j\leq w_j\quad \text{for}\;j=1,\ldots,m\,,\\
  v<w \quad &\Longleftrightarrow\quad\,v\,\leq\, w\quad\text{and} \quad
  v_j<w_j\quad\text{for some}\;j\in\{1,\ldots,m\}\,,
\end{split}
\]
We write $A\leq B$ for $m\times m$ matrices $A=[a_{ij}]_{ij}$ and
$B=[b_{ij}]_{ij}$ whenever $a_{ij}\leq b_{ij}$ for all $i,j$. Let $A$ be an
$m\times m$ quasipositive matrix, i.e. a matrix such that there exists
$\lambda>0$ with $A+\lambda I\geq 0$. Let $\rho >0$; let us recall the
definitions of exponential ordering on $BU$ given in~\cite{noov}. If
$x,\,y\in BU$, then
\[
\begin{split}
  x\leq_{A,\rho } y \;&\Longleftrightarrow\; x\leq y \;\text{ and }\;
  y(t)-x(t)\geq e^{A(t-s)}(y(s)-x(s))\,,
  -\rho \leq s\leq t\leq 0\,,\\
  x<_{A,\rho } y \;&\Longleftrightarrow\; x \leq_{A,\rho } y
  \quad\text{and}\quad x\neq
  y\,,\;\text{and} \\[.1cm]
  x\leq_{A,\infty} y\;&\Longleftrightarrow\; x\leq y \;\text{ and }\;
  y(t)-x(t)\geq e^{A(t-s)}(y(s)-x(s))\,,
  -\infty<s\leq t\leq 0\,,\\
  x<_{A,\infty} y \;&\Longleftrightarrow\; x \leq_{A,\infty}
  y\;\;\text{and}\quad x\neq y\,.
\end{split}
\]
In what follows, $\leq_A$ will denote any of the order relations
$\leq_{A,\rho}$ and $\leq_{A,\infty}$. However, in the case of
$\leq_{A,\infty}$, we will assume without further notice that all the
eigenvalues of $A$ have strictly negative real parts. The theory will
provide different dynamical conclusions for each choice. The aforementioned
relations define positive cones in $BU$, $BU^+_A=\{x\in BU:x\geq_A 0\}$, in
the sense that they are closed subsets of $BU$ and satisfy
$BU^+_A+BU^+_A\subset BU^+_A$, $\R^+BU^+_A\subset BU^+_A$ and
$BU^+_A\cap(-BU^+_A)=\{0\}$. Note that, if $\leq_A=\leq_{A,\rho}$, then a
smooth function (resp. a Lipschitz continuous function) $x$ belongs to
$BU_A^+$ if and only if $x\geq 0$ and $x'(s)\geq Ax(s)$ for each
(resp. a.e.) $s\in[-\rho,0]$, and, if $\leq_A=\leq_{A,\infty}$, then it
belongs to $BU_A^+$ if and only if $x\geq 0$ and $x'(s)\geq Ax(s)$ for each
(resp. a.e.) $s\in(-\infty,0]$.\par
On each fiber of the product $\Om\times BU$, we define the following
\emph{transformed exponential order} relation: if
$(\w,x),\,(\w,y)\in\Om\times BU$, then
\[
(\w,x)\leq_{D,A}(\w,y) \;\Longleftrightarrow\; \widehat
D_2(\w,x)\leq_A\widehat D_2(\w,y).
\]
Let us assume the following hypothesis:
\begin{itemize}
\item[(F1)] $F:\Om\times BU\to\R^m$ is continuous on $\Om\times BU$ and its
  restriction to $\Om\times B_r$ is Lipschitz continuous in its second
  variable when the norm is considered on $B_r$ for all $r>0$.
\end{itemize}
As seen in Wang and Wu~\cite{wawu} and Wu~\cite{wu}, for each $\w\in\Om$,
the local existence and uniqueness of the solutions of
equation~\ref{original_family} follow from (F1). Moreover, given
$(\w,x)\in\Om\times BU$, if $z({\cdot},\w,x)$ represents the solution of
equation~\ref{original_family} with initial datum $x$, then the mapping
$u(t,\w,x):(-\infty,0]\to\R^m$, $s\mapsto z(t+s,\w,x)$ is an element of
$BU$ for all $t\geq 0$ where $z({\cdot},\w,x)$ is defined.\par
Therefore, a local skew-product semiflow on $\Om\times BU$ can be defined
as follows:
\[
\begin{array}{rccl}
  \tau:&\mathcal U\subset\R^+\times\Om\times BU&\longrightarrow&\Om\times BU\\
  &(t,\w,x)&\mapsto&(\w{\cdot}t,u(t,\w,x)).
\end{array}
\]
Let $(\w,y)\in\Om\times BU$. For each $t\geq 0$ where $u(t,\wh
D^{-1}(\w,y))$ is defined, we define $\wh u(t,\w,y)=\wh
D_2(\w{\cdot}t,u(t,\wh D^{-1}(\w,y)))$. Let us check that
\[
\wh z({\cdot},\w,y):t\mapsto
\begin{cases}
  \,y(t)&\text{if }t\leq 0,\\
  \,\wh u(t,\w,y)(0)&\text{if }t\geq 0,
\end{cases}
\]
is the solution of \addtocounter{equation}{1}
\begin{equation}\tag*{(\arabic{section}.\arabic{equation})$_\w$}\label{transformed_family}
  \wh z\,'(t)=G(\w{\cdot}t,\wh z_t),\,t\geq 0,\,\w\in\Om
\end{equation}
through $(\w,y)$, where $G=F\circ\wh D^{-1}$. Let $x=(\wh D^{-1})_2(\w,y)$;
if $t\geq 0$, then
\[
\begin{split}
  \frac{d}{dt}\wh z(t,\w,y)&=\frac{d}{dt}\left[\wh
    D_2(\w{\cdot}t,u(t,\w,x))(0)\right]=\frac{d}{dt}D(\w{\cdot}t,u(t,\w,x))=F(\w{\cdot}t,u(t,\w,x))\\
  &=F\circ\wh D^{-1}(\wh D(\w{\cdot}t,u(t,\w,x)))=F\circ\wh
  D^{-1}(\w{\cdot}t,\wh u(t,\w,x)).
\end{split}
\]
It only remains to observe that a simple calculation yields $\wh
z({\cdot},\w,y)_t=\wh u(t,\w,y)$ for all $t\geq 0$. Let us assume some more
hypotheses concerning the map $F$.
\begin{itemize}
\item[(F2)] $F(\Om\times B_r)$ is a bounded subset of $\R^m$ for all $r>0$.
\item[(F3)] The restriction of $F$ to $\Om\times B_r$ is continuous when
  the compact-open topology is considered on $B_r$, for $r>0$.
\item[(F4)] If $(\w,x),\,(\w,y)\in\Om\times BU$ and
  $(\w,x)\leq_{D,A}(\w,y)$, then $F(\w,y)-F(\w,x)\geq A(D(\w,y)-D(\w,x))$.
\end{itemize}
\begin{proposition}\label{transformed_properties_1-4}
  Under hypotheses {\upshape (F1)--(F4)}, the following assertions hold:
  \begin{itemize}
  \item[(i)] $G$ is continuous on $\Om\times BU$ and its restriction to
    $\Om\times B_r$ is Lipschitz continuous in its second variable when the
    norm is considered on $B_r$ for all $r>0$;
  \item[(ii)] $G(\Om\times B_r)$ is a bounded subset of $\R^m$ for all
    $r>0$;
  \item[(iii)] the restriction of $G$ to $\Om\times B_r$ is continuous when
    the compact-open topology is considered on $B_r$, for $r>0$;
  \item[(iv)] if $(\w,x),\,(\w,y)\in\Om\times BU$ and $x\leq_Ay$, then
    $G(\w,y)-G(\w,x)\geq A(y(0)-x(0))$.
  \end{itemize}
\end{proposition}
\begin{proof}
  First, let us check (i). Let $\{(\w_n,x_n)\}_n\subset\Om\times BU$ be a
  sequence with $\w_n\to\w$ and $\n{x_n-x}_\infty\to 0$ as $n\to\infty$ for
  some $(\w,x)\in\Om\times BU$. Then $x_n\stackrel{\di}\to x$ as
  $n\to\infty$ and there is an $r>0$ such that $x_n,\,x\in B_r$ for all
  $n\in\N$. Consequently, $\wh D^{-1}(\w_n,x_n)\stackrel{\di}\to\wh
  D^{-1}(\w,x)$ as $n\to\infty$ and, from Theorem~\ref{D-1}, $\wh
  D^{-1}(\w_n,x_n),\,\wh D^{-1}(\w,x)\in B_{kr}$ for all $n\in\N$. Thanks
  to (F3), $G(\w_n,x_n)\to G(\w,x)$ as $n\to\infty$. As for the Lipschitz
  continuity, let $r>0$ and fix $(\w,y_1),\,(\w,y_2)\in\Om\times B_r$;
  denote by $x_i=(\wh D^{-1})_2(\w,y_i)$, $i=1,\,2$. From
  Theorem~\ref{D-1}, it follows that $\n{x_i}_\infty\leq
  k\n{y_i}_\infty\leq kr$, $i=1,\,2$. Let $L>0$ be the Lipschitz constant
  of $F$ on $\Om\times B_{kr}$. Again Theorem~\ref{D-1} yields
  \[
  \n{G(\w,y_1)-G(\w,y_2)}=\n{F(\w,x_1)-F(\w,x_2)}\leq
  L\n{x_1-x_2}_\infty\leq Lk\n{y_1-y_2}_\infty,
  \]
  and the result is proved.\\
  As for (ii), let $r>0$; from Theorem~\ref{D-1}, it follows that
  \[
  G(\Om\times B_r)=F(\wh D^{-1}(\Om\times B_r))\subset F(\Om\times B_{kr})
  \]
  and the latter is bounded thanks to (F2).\\
  Let us focus on (iii). Let $r>0$; once more, Theorem~\ref{D-1} implies
  that $\wh D^{-1}(\Om\times B_r)\subset\Om\times B_{kr}$, and $F$ is
  continuous
  there when we consider the compact-open topology.\\
  Finally, (iv) is a straightforward consequence of (F4).
\end{proof}
We may now define another local skew-product semiflow on $\Om\times BU$
from the solutions of the equations of the family~\ref{transformed_family}
(see Hino, Murakami and Naito~\cite{himn}) in the following manner:
\[
\begin{array}{rccl}
  \wh\tau:&\wh{\mathcal U}\subset\R^+\times\Om\times BU&\longrightarrow&\Om\times BU\\
  &(t,\w,x)&\mapsto&(\w{\cdot}t,\wh u(t,\w,x)).
\end{array}
\]
Given $r>0$, a forward orbit $\{\wh\tau(t,\w_0,x_0):t\geq 0\}$ of the
transformed skew-product semiflow $\wh\tau$ is said to be \emph{uniformly
  stable for the order $\leq_A$ in $B_r$} if, for every $\varepsilon>0$,
there is a $\delta>0$, called the \emph{modulus of uniform stability}, such
that, if $s\geq 0$ and $\di(\wh u(s,\w_0,x_0),x)\leq \delta$ for certain
$x\in B_r$ with $x\leq_A \wh u(s,\w_0,x_0)$ or $\wh u(s,\w_0,x_0)\leq_A x$,
then for each $t\geq 0$,
\[
\di(\wh u(t+s,\w_0,x_0),\wh u(t,\w_0{\cdot}s,x))= \di(\wh
u(t,\w_0{\cdot}s,\wh u(s,\w_0,x_0)),\wh
u(t,\w_0{\cdot}s,x))\leq\varepsilon.
\]
If this happens for each $r>0$, the forward orbit is said to be
\emph{uniformly stable for the order $\leq_A$ in bounded sets}.\par
Notice that an argument similar to the one given in Proposition~4.1
in~\cite{noos} and Proposition~4.2 in~\cite{muno} ensures that, for all
$(\w_0,x_0)\in\Om\times BU$ giving rise to a bounded solution,
$\cls_X\{u(t,\w_0,x_0):t\geq 0\}$ is a compact subset of $BU$ and the
omega-limit set of $(\w_0,x_0)$ can be defined as
\[
\mathcal O(\w_0,x_0)=\{(\w,x)\in\Om\times BU:\exists
t_n\uparrow\infty\text{ with
}\w_0{\cdot}t_n\to\w,\,u(t_n,\w_0,x_0)\stackrel{\di}\to x\}.
\]\par
Moreover, the restriction of $\tau$ to $\mathcal O(\w_0,x_0)$ is continuous
when the compact-open topology is considered on $BU$, and $\mathcal
O(\w_0,x_0)$ admits a flow extension. The main objective of this section is
to transfer the dynamical structure of the semiflow $(\Om\times
BU,\wh\tau,\R^+)$ to $(\Om\times BU,\tau,\R^+)$.
\begin{theorem}\label{monotone_tau}
  Fix $(\w,x),\,(\w,y)\in\Om\times BU$ such that
  $(\w,x)\leq_{D,A}(\w,y)$. Then
  \[
  \tau(t,\w,x)\leq_{D,A}\tau(t,\w,y)
  \]
  for all $t\geq 0$ where they are defined.
\end{theorem}
\begin{proof}
  It is clear that $\wh D_2(\w,x)\leq_A\wh D_2(\w,y)$. Now, from Theorem
  3.5 in~\cite{noov} and Proposition~\ref{transformed_properties_1-4}, it
  follows that $\wh u(t,\wh D(\w,x))\leq_A\wh u(t,\wh D(\w,y))$ or,
  equivalently, $\wh D_2(\w{\cdot}t,u(t,\w,x))\leq_A\wh
  D_2(\w{\cdot}t,u(t,\w,y))$ whenever they are defined. Therefore, we have
  $\tau(t,\w,x)\leq_{D,A}\tau(t,\w,y)$ for all $t\geq 0$ where they are
  defined, as wanted.
\end{proof}
Let us assume two more hypotheses concerning $F$ and the semiflow
$\wh\tau$. The fact that we are imposing a condition on the semiflow
$\wh\tau$ seems to suggest that such condition should be difficult to
verify when studying specific systems of equations. As it will be shown
later on, this kind of condition arises naturally in some systems and is
easier to check.
\begin{itemize}
\item[(F5)] There exists $r_0>0$ such that all the trajectories for
  $\wh\tau$ with a Lipschitz continuous initial datum within $B_{\wh r_0}$
  are relatively compact for the product metric topology and uniformly
  stable for the order $\leq_A$ in bounded sets, where
  \[
  \wh r_0=\n{A^{-1}}(\sup\{\n{F(\w,x)}:(\w,x)\in\wh D^{-1}(\Om\times
  B_{r_0})\}+\n{A}r_0).
  \]
\item[(F6)] If $(\w,x),\,(\w,y)\in\Om\times BU$ admit a backward orbit
  extension for the semiflow $\tau$, $(\w,x)\leq_{D,A}(\w,y)$ and there
  exists $J\subset\{1,\ldots,m\}$ such that
  \[
  \begin{split}
    &\phantom{xxx}\wh D_2(\w,x)_i=\wh D_2(\w,y)_i\,\text{ for all }i\notin
    J\,\text{
      and}\\
    &\phantom{xxx}\wh D_2(\w,x)_i(s)<\wh D_2(\w,y)_i(s)\,\text{ for all
    }i\in J\,\text{ and all }s\leq 0,
  \end{split}
  \]
  then $F_i(\w,y)-F_i(\w,x)-[A(D(\w,y)-D(\w,x))]_i>0$ for all $i\in J$.
\end{itemize}
The next result is an immediate consequence of these two properties, and so
we give it without a proof.
\begin{proposition}\label{transformed_properties_5-6}
  Under hypotheses {\upshape (F5)} and {\upshape (F6)}, the following
  assertions hold:
  \begin{itemize}
  \item[(i)] there exists $r_0>0$ such that all the trajectories for
    $\wh\tau$ with a Lipschitz continuous initial datum within $B_{\wh
      r_0}$ are relatively compact for the product metric topology and
    uniformly stable for the order $\leq_A$ in bounded sets, where
    \begin{equation}\label{radius_F5}
      \wh r_0=\n{A^{-1}}(\sup\{\n{G(\w,x)}:(\w,x)\in\Om\times
      B_{r_0}\}+\n{A}r_0).
    \end{equation}
  \item[(ii)] if $(\w,x),\,(\w,y)\in\Om\times BU$ admit a backward orbit
    extension for the semiflow $\wh\tau$, $x\leq_Ay$ and there exists
    $J\subset\{1,\ldots,m\}$ such that
    \[
    \begin{split}
      &x_i=y_i\,\text{ for all }i\notin J\,\text{ and}\\
      &x_i(s)<y_i(s)\,\text{ for all }i\in J\,\text{ and all }s\leq 0,
    \end{split}
    \]
    then $G_i(\w,y)-G_i(\w,x)-[A(y(0)-x(0))]_i>0$ for all $i\in J$.
  \end{itemize}
\end{proposition}
Relation~\eqref{radius_F5} is an improved version of formula (5.1)
in~\cite{noov}. Following this paper, we come now to the main result of
this section, which establishes the 1-covering property of the omega-limit
sets.
\begin{theorem}\label{copy_base}
  Assume that conditions {\upshape (D1)--(D3)} are satisfied and that $D$
  is stable; furthermore, assume conditions {\upshape (F1)--(F6)}.  Fix
  $(\w_0,x_0)\in\wh D^{-1}(\Om\times B_{r_0})$ such that $\{\wh\tau(t,\wh
  D(\w_0,x_0)):t\geq 0\}$ is relatively compact for the product metric
  topology and uniformly stable for $\leq_A$ in bounded sets, and such that
  $K=\mathcal O(\w_0,x_0)\subset\wh D^{-1}(\Om\times B_{r_0})$. If
  $\,\leq_A=\leq_{A,\infty}$, then we will further assume that $\wh
  D_2(\w_0,x_0)$ is Lipschitz continuous. Then $K=\{(\w,c(\w)):\w\in\Om\}$
  and
  \[
  \lim_{t\to\infty}\text{\sf d}(u(t,\w_0,x_0),c(\w_0{\cdot}t))=0,
  \]
  where $c:\Om\to BU$ is a continuous equilibrium,
  i.e. $c(\w{\cdot}t)=u(t,\w,c(\w))$ for all $t\geq 0$ and all $\w\in\Om$,
  considering on $BU$ the compact-open topology.
\end{theorem}
\begin{proof}
  As seen in Theorem 5.6 in~\cite{noov}, from
  Propositions~\ref{transformed_properties_1-4}
  and~\ref{transformed_properties_5-6}, it follows that $\mathcal O(\wh
  D(\w_0,x_0))=\{(\w,\wh c(\w)):\w\in\Om\}$ and
  \[
  \lim_{t\to\infty}\di(\wh u(t,\w_0,x_0),\wh c(\w_0{\cdot}t))=0,
  \]
  where $\wh c:\Om\to BU$ is a continuous equilibrium considering on $BU$
  the compact-open topology. We observe that
  \[
  K\!=\!\wh D^{-1}(\mathcal O(\wh D(\w_0,x_0)))\!=\!\wh D^{-1}(\{(\w,\!\wh
  c(\w))\!:\!\w\in\Om\})=\{(\w,(\wh D^{-1})_2(\w,\!\wh
  c(\w)))\!:\!\w\in\Om\}.
  \]
  Let us define $c:\Om\to BU$, $\w\mapsto(\wh D^{-1})_2(\w,\wh c(\w))$. The
  continuity of $c$ when we consider the compact-open topology on $BU$ is a
  consequence of Theorem~\ref{D-1} and the fact that $\mathcal O(\wh
  D(\w_0,x_0))\subset\Om\times B_{r_0}$. Moreover, $c$ is an equilibrium,
  for its graph defines an omega-limit set. Finally, from Theorem~\ref{D-1}
  again and the boundedness of the trajectory and of $\wh c(\Om)$, we
  conclude that
  \[
  \lim_{t\to\infty}\di(u(t,\w_0,x_0),c(\w_0{\cdot}t))=0,
  \]
  and the proof is finished.
\end{proof}
\section{Neutral compartmental systems}\label{section_comp_sys}
Let $(\Om,d)$ be a compact metric space and let $\R\times\Om\to\Om$,
$\w\mapsto\w{\cdot}t$ be a minimal real flow on $\Om$. In this section, we
apply the foregoing results to the study of compartmental models, used to
describe processes in which the transport of material among some
compartments takes a non-negligible length of time, and each compartment
produces or swallows material.\par
Let us suppose that we have a system formed by $m$ compartments
$C_1,\ldots,C_m$, and denote by $z_i(t)$ the amount of material within
compartment $C_i$ at time $t$ for each $i\in\{1,\ldots,m\}$. Material flows
from compartment $C_j$ into compartment $C_i$ through a pipe having a
transit time distribution given by a positive regular Borel measure
$\mu_{ij}$ with total variation $\mu_{ij}((-\infty,0])=1$, for each $i$,
$j\in\{1,\ldots,m\}$, whereas the outcome of material from $C_i$ to $C_j$
is assumed to be instantaneous. Let
$g_{ij}:\Om\times\mathbb{R}\to\mathbb{R}$ be the so-called \emph{transport
  function} determining the volume of material flowing from $C_j$ to $C_i$
given in terms of the time and the amount of material within $C_j$ for
$i,\,j\in\{1,\ldots,m\}$. For each $i\in\{1,\ldots,m\}$, there is an inflow
of material from the environment given by $I_i:\Om\to\R$, an outflow of
material toward the environment given by
$g_{0i}:\Om\times\mathbb{R}\to\mathbb{R}$, and the compartment $C_i$
produces or swallows material itself at a rate given by some regular Borel
measures $\nu_{ij}(\w)$, $\w\in\Om$, and some functions $b_{ij}:\Om\to\R$,
$j\in\{1,\ldots,m\}$.\par
Once the destruction and creation of material is taken into account, the
change of the amount of material of any compartment $C_i$, $1\leq i\leq m$,
equals the difference between the amount of total inflow into and total
outflow out of $C_i$, and we obtain a model governed by the following
system of NFDEs: \addtocounter{equation}{1}
\begin{equation}\tag*{(\arabic{section}.\arabic{equation})$_\w$}\label{compartmental_system}
  \begin{split}
    \frac{d}{dt}&\sum_{j=1}^m\left[b_{ij}(\w{\cdot}t)z_j(t)-\int_{-\infty}^0
      z_j(t+s)\,d\nu_{ij}(\w{\cdot}t)(s)\right]=-g_{0i}(\w{\cdot}t,z_i(t))\\
    -&\sum_{j=1}^m g_{ji}(\w{\cdot}t,z_i(t)) + \sum_{j=1}^m\int_{-\infty}^0
    g_{ij}(\w{\cdot}(t+s),z_j(t+s))\,d\mu_{ij}(s)+I_i(\w{\cdot}t),
  \end{split}
\end{equation}
for $t\geq 0$, $i=1,\ldots,m$ and $\w\in\Om$, where
$g_{ij}:\Om\times\R\to\R$, $I_i:\Om\to\R$ and $\mu_{ij}$, $\nu_{ij}(\w)$
are regular Borel measures with finite total variation for all
$\w\in\Om$. For the sake of simplicity, let us denote
$g_{i0}:\Om\times\R\to\R$, $(\w,v)\mapsto I_i(\w)$ and let
$g=(g_{ij})_{ij}:\Om\times\R\to \R^{m^2+2m}$. We denote by $B(\w)$ and
$\nu(\w)$ the matrices $[b_{ij}(\w)]_{ij}$ and $[\nu_{ij}(\w)]_{ij}$,
$\w\in\Om$, respectively.\par
Let $F:\Om\times BU\to\R^m$ be the map defined for $(\w,x)\in\Om\times BU$
by
\[
F_i(\w,x)=-\sum_{j=0}^m g_{ji}(\w,x_i(0)) + \sum_{j=1}^m\int_{-\infty}^0
g_{ij}(\w{\cdot}s,x_j(s))\,d\mu_{ij}(s)+I_i(\w),
\]
$i=1,\ldots,m$, and let $D:\Om\times BU\to\R^m$ be the map defined for
$(\w,x)\in\Om\times BU$ by
\[
D(\w,x)=\left(\sum_{j=1}^m\left[b_{ij}(\w)x_j(0)-\int_{-\infty}^0x_j\,d\nu_{ij}(\w)\right]\right)_{i=1}^m
=B(\w)x(0)-\int_{-\infty}^0[d\nu(\w)]x.
\]
With this notation, the family of equations~\ref{compartmental_system} can
be written as \addtocounter{equation}{1}
\begin{equation}\tag*{(\arabic{section}.\arabic{equation})$_\w$}\label{compartmental_short}
  \frac{d}{dt}D(\w{\cdot}t,z_t)=F(\w{\cdot}t,z_t),\quad t\geq 0,\,\w\in\Om.
\end{equation}
We will assume the following hypotheses on the
family~\ref{compartmental_system}.
\begin{itemize}
\item[(C1)] $g_{ij}$ is $C^1$ and nondecreasing in its second variable, and
  $g_{ij}(\w,0)=0$ for each $\w\in\Om$, $i\in\{0,\ldots,m\}$ and
  $j\in\{1,\ldots,m\}$.
\item[(C2)] $\mu_{ij}((-\infty,0])=1$ and $\int_{-\infty}^0
  |s|\,d\mu_{ij}(s)<\infty$ for $i,j=1,\ldots,m$.
\item[(C3)] $\nu_{ij}(\w)(\{0\})=0$ for all $\w\in\Om$ and the mapping
  $\nu:\Om\to\mathcal{M}$, $\w\mapsto\nu(\w)$ is continuous, where
  $\mathcal{M}$ is the Banach space of $m\times m$ matrices of Borel
  measures on $(-\infty,0]$ with the supremum norm defined from the total
  variation of the measures.
\item[(C4)] $B(\w)$ is a regular matrix for all $\w\in\Om$ and $B:\Om\to
  \mathbb M_m(\R)$, $\w\mapsto B(\w)$ is continuous; moreover,
  $B(\w)^{-1}\geq 0$, $B(\w)^{-1}\nu(\w)$ is a matrix of positive measures
  for all $\w\in\Om$, and
  \[
  \n{B(\w)^{-1}\nu(\w)}_\infty((-\infty,0])<1.
  \]
\end{itemize}\par
First, in this section we will prove the stability of the operator
$D$. Later, we will deduce that at least the trajectories transformed by
the operator $\wh D$ with Lipschitz continuous initial data are uniformly
stable for $\leq_A$ on bounded sets.\par
Conditions (C1)--(C4) yield the following result.
\begin{proposition}
  The mapping $D$ satisfies {\upshape (D1)--(D3)} and the mapping $F$
  satisfies {\upshape (F1)--(F3)} respectively.
\end{proposition}
We define the map $\wh D:\Om\times BU\to\Om\times BU$ as in
Theorem~\ref{Dhat_cont}. For each $\w\in\Om$, let us define $\wh B_\w:BU\to
BU$, $\wh B_\w(x)(s)=B(\w{\cdot}s)x(s)$. It is easy to check that $\wh
B_\w$ is a linear isomorphism of $BU$ and it is continuous for the norm; in
addition $(\w,x)\mapsto\wh B_\w(x)$ is continuous on each set of the form
$\Om\times B_r$,
$r>0$, when the compact-open topology is considered. We can denote $(\wh
B_\w)^{-1}=(\wh{B^{-1}})_\w$, $\w\in\Om$.
\begin{theorem}\label{explicit_D-1}
  For each $\w\in\Om$, consider the continuous linear operator $\wh
  L_\w:BU\to BU$ defined for $x\in BU$ and $s\leq 0$ by
  \[
  \wh L_\w(x)(s)=B(\w{\cdot}s)^{-1}\int_{-\infty}^0[d\nu(\w{\cdot}s)]x_s.
  \]
  Then the following statements hold:
  \begin{itemize}
  \item[(i)] $\sup_{\w\in\Om}\n{\wh L_\w}<1$ and $\wh D_2(\w,x)=[\wh
    B_\w\circ(I-\wh L_\w)](x)$ for each $(\w,x)\in\Om\times BU$;
  \item[(ii)] $\wh D$ is invertible and \vspace{-2mm}
    \[
    (\wh D^{-1})_2(\w,x)=\sum_{n=0}^\infty(\wh L_\w^n\circ(\wh
    B_\w)^{-1})(x) \vspace{-1mm}
    \]
    for every $(\w,x)\in\Om\times BU$;
  \item[(iii)] $(\wh D^{-1})_2(\w,x)\geq 0$ for every $(\w,x)\in\Om\times
    BU$ with $x\geq 0$;
  \item[(iv)] the map $\Om\times B_r\to B_r$, $(\w,x)\mapsto\wh L_\w(x)$ is
    uniformly continuous for the compact-open topology for each $r>0$;
  \item[(v)] $D$ is stable.
  \end{itemize}
\end{theorem}
\begin{proof}
  From condition (C4), we conclude that $\sup_{\w\in\Om}\n{\wh L_\w}<1$. It
  is immediate to check that $\wh D_2(\w,x)=[\wh B_\w\circ(I-\wh L_\w)](x)$
  for each $(\w,x)\in\Om\times BU$. This way, (i) is proved, whence $\wh D$
  is invertible and \vspace{-2mm}
  \[
  (\wh D^{-1})_2(\w,x)=\sum_{n=0}^\infty(\wh L_\w^n\circ(\wh B_\w)^{-1})(x)
  \vspace{-1mm}
  \]
  for every $(\w,x)\in\Om\times BU$, which proves (ii). Now, (iii) follows
  from (ii) and hypothesis (C4).\par
  In addition, for all $r_1>0$ and all $(\w,x)\in\Om\times B_{r_1}$, it is
  clear that $(\wh B_\w)^{-1}(x)\in\Om\times B_{r_2}$ and $\wh
  L_\w(x)\in\Om\times B_{r_1}$, where
  $r_2=\sup_{\w_1\in\Om}\n{B(\w_1)^{-1}}r_1$; besides, if $(\w_1,x_1)$,
  $(\w_2,x_2)\in\Om\times B_{r_1}$ and $s\in[-\rho,0]$ for some $\rho>0$,
  then
  \[
  \begin{split}
    \n{\wh L_{\w_1}&(x_1)(s)\!-\!\wh L_{\w_2}(x_2)(s)}\!\leq\!\n{\wh
      L_{\w_1}(x_1)(s)\!-\!\wh L_{\w_1}(x_2)(s)}\!+\!\n{\wh
      L_{\w_1}(x_2)(s)-\wh L_{\w_2}(x_2)(s)}\\
    \leq&\sup_{\w\in\Om}\n{B(\w)^{-1}}\sup_{\w\in\Om}\n{\nu(\w)}_\infty((-\infty,-\rho])2r_1+\n{x_1-x_2}_{[-2\rho,0]}\\
    &+\n{B(\w_1{\cdot}s)^{-1}-B(\w_2{\cdot}s)^{-1}}\sup_{\w\in\Om}\n{\nu(\w)}_\infty((-\infty,0])r_1\\
    &+\sup_{\w\in\Om}\n{B(\w)^{-1}}\n{\nu(\w_1{\cdot}s)-\nu(\w_2{\cdot}s)}_\infty((-\infty,0])r_1.
  \end{split}
  \]
  This proves (iv). Now it is immediate to deduce that, given $r>0$, $\wh
  D^{-1}$ is uniformly continuous on $\Om\times B_r$ for the compact-open
  topology. Let us fix $\eps>0$ and $\rho>0$. There is an $n_0\in\N$ such
  that $\sum_{n=n_0}^\infty\n{\wh L_\w^n\circ(\wh
    B_\w)^{-1}}<\eps/(3r)$. From (iv) and the continuity of
  $(\w,x)\mapsto\wh B_\w(x)$ for the product metric topology on each ball,
  it follows that there exist $\rho_0>0$ and $\delta>0$ such that, if
  $(\w_1,x_1)$, $(\w_2,x_2)\in\Om\times B_r$ with $d(\w_1,\w_2)<\delta$ and
  $\n{x_1(s)-x_2(s)}<\delta$ for all $s\in[-\rho_0,0]$, then $\n{\wh
    L_{\w_1}^j\circ(\wh B_{\w_1})^{-1}(x_1)(s)-\wh L_{\w_2}^j\circ(\wh
    B_{\w_2})^{-1}(x_2)(s)}\leq\eps/(3n_0)$, $j\in\{0,\ldots,n_0-1\}$,
  $s\in[-\rho,0]$. As a consequence, $\left\|\sum_{n=0}^\infty\wh
    L_{\w_1}^n\circ(\wh B_{\w_1})^{-1}(x_1)(s)-\sum_{n=0}^\infty\wh
    L_{\w_2}^n\circ(\wh B_{\w_2})^{-1}(x_2)(s)\right\|\leq\eps$ for every
  $s\in[-\rho,0]$, which proves, according to Theorem~\ref{D-1}, that $D$
  is stable. This completes the proof.
\end{proof}
We will assume now that (F4) is satisfied. Some sufficient conditions for
(F4) to hold will be studied in the next section.\par
Consider the \emph{total mass} of the system~\ref{compartmental_system},
$M:\Om\times BU\to\R$, defined for $(\w,x)\in\Om\times BU$ by
\[
M(\w,x)=\sum_{i=1}^mD_i(\w,x)+\sum_{i=1}^m\sum_{j=1}^m \int_{-\infty}^0
\left( \int_s^0g_{ji}(\w{\cdot} \tau,x_i(\tau))d\tau \right)d\mu_{ji}(s).
\]
$M$ is well defined due to (C2) and because, if $i,j\in\{1,\ldots,m\}$ and
$(\w,x)\in\Om\times BU$, then
$\left|\int_s^0g_{ji}(\w{\cdot}\tau,x_i(\tau))d\tau\right|\leq c_1|s|$,
where $c_1$ is a bound of $g_{ji}$ on
$\Om\times[-\n{x}_\infty,\n{x}_\infty]$.\par
The next result is in the line of some results found in~\cite{muno},
\cite{noov} and \cite{wufr}, and states an important equality satisfied by
the total mass of a compartmental system.
\begin{proposition}\label{M_variation}
  $M$ is uniformly continuous on $\Om\times B_r$ for all $r>0$ for the
  product metric topology. Moreover,
  \[
  M(\tau(t,\w,x))=M(\w,x)+\sum_{i=1}^m\int_0^t(I_i(\w{\cdot}s)-g_{0i}(\w{\cdot}s,z_i(s,\w,x)))\,ds
  \]
  for all $t\geq 0$ where the solution is defined and all
  $(\w,x)\in\Om\times BU$.
\end{proposition}
\begin{lemma}\label{M_bounds_D}
  Fix $(\w,x)$, $(\w,y)\in\Om\times BU$ with $(\w,x)\leq_{D,A}(\w,y)$. Then
  \[
  0\leq D_i(\tau(t,\w,y))-D_i(\tau(t,\w,x))\leq M(\w,y)-M(\w,x)
  \]
  for each $i=1,\ldots,m$ and whenever $z(t,\w,x)$ and $z(t,\w,y)$ are
  defined.
\end{lemma}
\begin{proof}
  It follows from (F4) and Theorem~\ref{monotone_tau} that the skew-product
  semiflow induced by~\ref{compartmental_system} is monotone. Hence, if
  $(\w,x)\leq_{D,A}(\w,y)$ then $\tau(t,\w,x)\leq_{D,A}\tau(t,\w,y)$
  whenever they are defined. From this and Theorem~\ref{explicit_D-1}, we
  also deduce that $x\leq y$ and $u(t,\w,x)\leq u(t,\w,y)$. Therefore,
  $D_i(\tau(t,\w,x))\leq D_i(\tau(t,\w,y))$ and $z_i(t,\w,x)\leq
  z_i(t,\w,y)$ for each $i=1,\ldots,m$. In addition, the monotonicity of
  transport functions yields $g_{ij}(\w{\cdot}t,z_j(t,\w,x))\leq
  g_{ij}(\w{\cdot}t,z_j(t,\w,y))$. From all these inequalities, the
  definition of total mass and Proposition~\ref{M_variation}, we deduce
  that
  \[
  \begin{split}
    0&\leq D_i(\tau(t,\w,y))-D_i(\tau(t,\w,x))\leq \sum_{i=1}^m
    \left[D_i(\tau(t,\w,y)) -D_i(\tau(t,\w,x))\right]\\
    &\leq M(\tau(t,\w,y))- M(\tau(t,\w,x))\\
    &=M(\w,y)-M(\w,x)+\sum_{i=1}^m\int_0^t
    (g_{0i}(\w{\cdot}s,z_i(s,\w,x))-g_{0i}(\w{\cdot}s,z_i(s,\w,y)))\,ds\\
    &\leq M(\w,y)-M(\w,x)\,,
  \end{split}
  \]
  as stated.
\end{proof}
\begin{proposition}\label{stability}
  Fix $r>0$. Given $\varepsilon>0$ there exists $\delta>0$ such that, if
  $(\w,x)$, $(\w,y)\in\Om\times B_r$ with {\upshape$\di(x,y)<\delta$} and
  $x\leq_Ay$, then $\|\wh z(t,\w,x)-\wh z(t,\w,y)\|\leq\varepsilon$
  whenever they are defined.
\end{proposition}
\begin{proof}
  Let $r_1=r\sup_{\w\in\Om}\n{(\wh D^{-1})_2(\w,{\cdot})}$. Fix $\eps>0$;
  it follows from Proposition~\ref{M_variation} that there is a
  $\delta_1>0$ such that, if $(\w,x)$, $(\w,y)\in\Om\times B_{r_1}$ with
  $\di(x,y)<\delta_1$, then $|M(\w,y)-M(\w,x)|\leq\eps$. Now, thanks to
  Theorem~\ref{explicit_D-1}, there is a $\delta>0$ such that, if $(\w,x)$,
  $(\w,y)\in\Om\times B_r$ with $\di(x,y)<\delta$, then
  \[
  \di((\wh D^{-1})_2(\w,x),(\wh D^{-1})_2(\w,y))<\delta_1.
  \]
  Altogether, using Lemma~\ref{M_bounds_D}, if $(\w,x)$,
  $(\w,y)\in\Om\times B_r$ with $\di(x,y)<\delta$ and $x\leq_Ay$, then
  \[
  0\leq D_i(\tau(t,\wh D^{-1}(\w,y)))-D_i(\tau(t,\wh D^{-1}(\w,x)))\leq
  M(\wh D^{-1}(\w,y))-M(\wh D^{-1}(\w,x))\leq\eps,
  \]
  whence
  \[
  0\leq\wh z_i(t,\w,y)-\wh z_i(t,\w,x)\leq\eps
  \]
  for all $t\geq 0$ where they are defined and all
  $i\in\{1,\ldots,m\}$. The result is proved.
\end{proof}
\begin{proposition}\label{1bounded_allbounded}
  Assume hypotheses {\upshape (C1)--(C4)} together with {\upshape (F4)}.
  Suppose that there exists $(\w_0,x_0)\in\Om\times BU$ such that
  $\wh\tau({\cdot},\w_0,y_0)$ is bounded, where $y_0=\wh
  D_2(\w_0,x_0)$. The following statements hold:
  \begin{itemize}
  \item[(i)] when the order $\leq_A$ associated to $\leq_{A,\rho}$ is
    considered, then it holds that, for all $(\w,x)\in\Om\times BU$,
    $\wh\tau({\cdot},\w,y)$ is bounded and uniformly stable for $\leq_A$ in
    bounded subsets, where $y=\wh D_2(\w,x)$;
  \item[(ii)] when the order $\leq_A$ associated to $\leq_{A,\infty}$ is
    considered, if $y_0$ is Lipschitz continuous, then it holds that, for
    all $(\w,x)\in\Om\times BU$ such that $y=\wh D_2(\w,x)$ is Lipschitz
    continuous, $\wh\tau({\cdot},\w,y)$ is bounded and uniformly stable for
    $\leq_A$ in bounded subsets.
  \end{itemize}
\end{proposition}
\begin{proof}
  To prove (i), let us define $\wt y\in BU$ as the solution of
  \begin{equation}\label{ord_eq}
    \left\{\begin{array}{ll}
        y'(t)=Ay(t)+\mathbf 1,&t\in[-\rho,0],\\
        y_{-\rho}\equiv\mathbf 1,
      \end{array}\right.
  \end{equation}
  where $\mathbf 1=(1,\ldots,1)\in\R^m$. Since the matrix $A$ is
  quasipositive, the system of ordinary differential
  equations~\eqref{ord_eq} is cooperative. The standard comparison theory
  for its solutions allows us to conclude that there exists $k_0>0$ such
  that $\wt y(t)\geq k_0\mathbf 1$ for all $t\leq 0$ and $\wt
  y\geq_{A,\rho}0$. Let us fix $z\in BU$ and suppose that $z$ is Lipschitz
  continuous on $[-\rho,0]$; let us check that there is a $\lambda_0>0$
  such that, for all $\lambda\geq\lambda_0$, $-\lambda\wt
  y\leq_{A,\rho}z\leq_{A,\rho}\lambda\wt y$. Since $\wt y(t)\geq k_0\mathbf
  1$ for all $t\leq 0$, it is clear that, if $\lambda\geq\lambda_0$, then
  $-\lambda\wt y\leq z\leq\lambda\wt y$ holds for all big enough
  $\lambda_0$. In addition, $(\lambda\wt y-z)'-A(\lambda\wt
  y-z)=\lambda\mathbf 1-(z'-Az)$ is greater or
  equal to 0 a.e. in $[-\rho,0]$ for all big enough $\lambda_0$ and we are done.\\
  Now, let $K$ be the omega-limit set of $(\w_0,y_0)$ for the semiflow
  $\wh\tau$ and let $r_1>0$ be such that $K\subset\Om\times B_{r_1}$. We
  will prove that $\wh z({\cdot},\w,\lambda\wt y)$ and $\wh
  z({\cdot},\w,-\lambda\wt y)$ are bounded for all $\w\in\Om$ and all
  sufficiently big $\lambda$. In order to do this, fix $\w\in\Om$ and $z\in
  K_\w$. We know that there is a $\lambda_0>0$ (irrespective of $\w$) such
  that, for all $\lambda\geq\lambda_0$, $-\lambda\wt
  y\leq_{A,\rho}z\leq_{A,\rho}\lambda\wt y$. For each $s\in[0,1]$, let
  $y_s=(1-s)z+s\lambda\wt y$; clearly, $y_s\leq_{A,\rho}y_t$ for all $0\leq
  s\leq t\leq 1$. Besides, there exists $r>0$ such that
  $\{y_s\}_{s\in[0,1]}\subset B_r$. An application of
  Proposition~\ref{stability} for $\eps=1$, implies that there are a
  $\delta>0$ and a partition $0=s_0\leq s_1\leq\cdots\leq s_n=1$ of $[0,1]$
  such that $\di(y_{s_j},y_{s_{j+1}})<\delta$ for all
  $j\in\{1,\ldots,n-1\}$, and therefore $\n{\wh z(t,\w,y_{s_j})-\wh
    z(t,\w,y_{s_{j+1}})}\leq 1$ for all $t\geq 0$ where they are
  defined. Consequently, for each $j\in\{0,\ldots,n\}$, the solution $\wh
  z({\cdot},\w,y_{s_j})$ is globally defined and $\n{\wh z(t,\w,z)-\wh
    z(t,\w,\lambda\wt y)}\leq n$ for all $t\geq 0$, which implies that $\wh
  z({\cdot},\w,\lambda\wt y)$ is bounded. Analogously, $\wh
  z({\cdot},\w,-\lambda\wt y)$ is bounded as well.\\
  Finally, let $(\w,x)\in\Om\times BU$ and $y=\wh D_2(\w,x)$. Let $z=\wh
  u(\rho,\w,y)$; since $z$ is Lipschitz continuous on $[-\rho,0]$, there
  exists $\lambda\geq\lambda_0$ such that $-\lambda\wt
  y\leq_{A,\rho}z\leq_{A,\rho}\lambda\wt y$, which implies that $\wh
  z({\cdot},\w{\cdot}\rho,z)$ is bounded thanks to the fact that $\wh
  z({\cdot},\w,\lambda\wt y)$ and $\wh z({\cdot},\w,-\lambda\wt y)$ are
  bounded and to the monotonicity of $\wh\tau$. Consequently, the
  trajectory through $(\w,y)$ for $\wh\tau$ is bounded. The remainder of
  the proof follows from Proposition~\ref{stability}.\par
  Now we deal with statement (ii). Notice that the fact that all the
  eigenvalues of $A$ have a negative real part implies that $A$ is a
  hyperbolic matrix. Consider the cooperative system of ordinary
  differential equations $y'=Ay+\mathbf 1$. It is well-known that there
  exists a unique solution of the aforementioned system that is bounded and
  exponentially stable when $t\to\infty$, namely $\wt y\equiv
  -A^{-1}\mathbf 1$. Since 0 is a strong subequilibrium of the system, as
  seen in Novo, N{\'u}{\~n}ez and Obaya~\cite{nono}, there exists $k_0>0$
  such that $\wt y\geq k_0\mathbf 1$. Denote again by $\wt y$ its
  restriction to $(-\infty,0]$. The rest of the proof is analogous to the
  one of statement (i).
\end{proof}
This proposition proves condition (F5) stated in
Section~\ref{section_copy_base}. The following result is a direct
consequence of Theorem~\ref{copy_base}.
\begin{theorem}\label{copy_base_comp}
  Under the hypotheses of {\upshape Proposition~\ref{1bounded_allbounded}},
  for all $(\w,x)\in\Om\times BU$ satisfying the conditions of statements
  {\upshape (i)} or {\upshape (ii)} respectively, we have that $\mathcal
  O(\w,x)$ is a copy of the base whenever {\upshape (F6)} holds as well.
\end{theorem}
Sufficient conditions under which (F6) holds will be given in the next
section together with those guaranteeing (F4).
\section{Some specific neutral compartmental
  systems}\label{section_examples}
Let $(\Om,d)$ be a compact metric space and let $\R\times\Om\to\Om$,
$\w\mapsto\w{\cdot}t$ be a minimal real flow on $\Om$. Let us focus on the
study of the following family of compartmental systems:
\addtocounter{equation}{1}
\begin{equation}\tag*{(\arabic{section}.\arabic{equation})$_\w$}\label{comp_diag}
  \frac{d}{dt}[z_i(t)-c_i(\w{\cdot}t)z_i(t\!-\alpha_i)]=\!-\!\!\sum_{j=1}^mg_{ji}(\w{\cdot}t,z_i(t))+\!\sum_{j=1}^mg_{ij}(\w{\cdot}(t\!-\!\rho_{ij}),z_j(t\!-\!\rho_{ij}))
\end{equation}
for $i\in\{1,\ldots,m\}$, $t\geq 0$ and $\w\in\Om$, where
$g_{ij}:\Om\times\R\to\R$ and $c_i:\Om\to\R$ are continuous functions and
$\alpha_i,\,\rho_{ij}\in\R$ for $i,\,j\in\{1,\ldots,m\}$. Throughout this
section, we will use the notation introduced in
Section~\ref{section_comp_sys}. Let us assume the following conditions on
equation~\ref{comp_diag}:
\begin{itemize}
\item[(G1)] $g_{ij}$ is $C^1$ and nondecreasing in its second variable, and
  $g_{ij}(\w,0)=0$ for each $\w\in\Om$, $i,\,j\in\{1,\ldots,m\}$;
\item[(G2)] $\alpha_i>0$, $\rho_{ij}\geq 0$ and $0\leq c_i(\w)<1$ for all
  $i,j\in\{1,\ldots,m\}$ and all $\w\in\Om$.
\end{itemize}
Notice that the family of equations~\ref{comp_diag} corresponds to a
\emph{closed} compartmental system, that is, a system where there is no
incoming material from the environment and there is no outgoing material
toward the environment either. As a result, the total mass of the system is
invariant along the trajectories and 0 is a constant bounded solution of
all the equations of the family.
\begin{definition}
  Let us consider a function $\mathfrak c:\Om\to \R^m$.
  \begin{itemize}
  \item[(i)] $\mathfrak c$ is said to be \emph{Lipschitz continuous along
      the flow} if, for some $\w\in\Om$, the mapping $\R\to\R^m$, $t\mapsto
    \mathfrak c(\w{\cdot}t)$ is Lipschitz continuous;
  \item[(ii)] $\mathfrak c$ is \emph{continuously differentiable along the
      flow} if, for all $\w\in\Om$, the mapping
    \[
    \begin{array}{rcl}
      \Om&\longrightarrow&\R^m\\
      \w&\mapsto&\left.\frac{d}{dt}\mathfrak c(\w{\cdot}t)\right|_{t=0}
    \end{array}
    \]
    is well-defined and continuous. We will refer to this mapping as the
    derivative of $\mathfrak c$.
  \end{itemize}
\end{definition}
It is easy to check that, due to the density of all trajectories within
$\Om$, if $\mathfrak c:\Om\to \R^m$ is Lipschitz continuous along the flow,
then the mappings $\R\to\R^m$, $t\mapsto \mathfrak c(\wt\w{\cdot}t)$,
$\wt\w\in\Om$, are all Lipschitz continuous with the same Lipschitz
constant.  From now on, $c:\Om\to\R^m$ will denote
$c=(c_i)_{i=1}^m:\w\mapsto(c_i(\w))_{i=1}^m$. It is noteworthy that,
given $\w\in\Om$, if $c$ is a Lipschitz continuous function, then $x\in BU$
is Lipschitz continuous if and only if $\wh D_2(\w,x)$ is Lipschitz
continuous.\par
For each $i\in\{1,\ldots,m\}$, let $c_i^{[0]}(\w)=1$ and, for each
$n\in\N$, define
\[
c_i^{[n]}(\w)=\prod_{j=0}^{n-1}c_i(\w{\cdot}(-j\alpha_i)),\,\w\in\Om.
\]
\begin{proposition}
  Assume that $c$ is continuously differentiable along the flow. Suppose
  that $(\w,x)\in\Om\times BU$ admits a backward orbit extension and that
  there is an $r_1>0$ such that $u(t,\w,x)\in B_{r_1}$ for each
  $t\in\R$. Then $z=z({\cdot},\w,x)$, the solution of~{\rm\ref{comp_diag}}
  with initial value $x$, belongs to $C^1(\R,\R^m)$.
\end{proposition}
\begin{proof}
  Let $\wh z=\wh z({\cdot},\wh D(\w,x))$. It is clear that $\wh z$ is of
  class $C^1$ and it is bounded by
  $r_1\sup_{\w_1\in\Om}\n{D(\w_1,{\cdot})}$. Then, for all $t\in\R$ and all
  $i\in\{1,\ldots,m\}$, it holds that
  \[
  z_i(t)=\sum_{n=0}^\infty c_i^{[n]}(\w{\cdot}t)\wh z_i(t-n\alpha_i).
  \]
  From (G2), it follows that this series converges uniformly on
  $\R$. Analogously, the formal derivative of the former series, namely
  \begin{equation}\label{deriv_omega_limit}
    \sum_{n=0}^\infty c_i^{[n]}(\w{\cdot}t)\wh
    z_i\,'(t-n\alpha_i)+\sum_{n=0}^\infty\left.\frac{d}{ds}c_i^{[n]}(\w{\cdot}(t+s))\right|_{s=0}\wh
    z_i(t-n\alpha_i),
  \end{equation}
  converges uniformly on $\R$ thanks to (G2). Consequently, $z_i$ is
  continuously differentiable on $\R$.
\end{proof}
Note that, in the conditions of the previous proposition, the derivative of
$z$ is given by~\eqref{deriv_omega_limit}. The following result is a
straightforward consequence of hypotheses (G1)--(G2).
\begin{proposition}\label{G1-G2_imply_C1-C4}
  Under hypotheses {\upshape(G1)--(G2)}, the family of
  equations~{\upshape\ref{comp_diag}} satisfies conditions
  {\upshape(C1)--(C4)}.
\end{proposition}
We next analyze some situations where the previous theory can be
applied. They are chosen to describe different types of conditions which
assure the monotonicity of the semiflow for some transformed exponential
order leading to different dynamical implications. In the next statement,
we take $\rho_{ii}=2\alpha_i$, $i=1,\ldots,m$. Note that condition (G3.2),
independently of the matrix $A$, is always required. The conclusion of the
theorem also guarantees that all the trajectories are relatively
compact.\par
For each $\w\in\Om$ and each $i,\,j\in\{1,\ldots,m\}$, let
$l_{ii}^-(\w)=\inf_{v\in\R}\frac{\partial g_{ii}}{\partial v}(\w,v)$,
$l_{ij}^+(\w)=\sup_{v\in\R}\frac{\partial g_{ij}}{\partial v}(\w,v)$ and
$L_i^+(\w)=\sum_{j=1}^ml_{ji}^+(\w)$. In this section, we will assume that
$L_i^+(\w)<\infty$ for all $\w\in\Om$.
\begin{theorem}\label{trans_rho=2alpha}
  Let us assume hypotheses {\upshape(G1)--(G2)} together with
  \begin{itemize}
  \item[(G3)] for each $i\in\{1,\ldots,m\}$, if $c_i\not\equiv 0$, then
    $\rho_{ii}=2\alpha_i$ and there exists $a_i\in(-\infty,0]$ such that,
    for all $\w\in\Om$, the following conditions hold:
    \begin{itemize}
    \item[(G3.1)] $(-a_i-L_i^+(\w))e^{a_i\alpha_i}-L_i^+(\w)c_i(\w)\geq 0$,
    \item[(G3.2)]
      $l_{ii}^-(\w{\cdot}(-\rho_{ii}))-L_i^+(\w)c_i^{[2]}(\w)\geq 0$,
    \end{itemize}
  \end{itemize}
  where at least one of the inequalities is strict. Then all the
  trajectories of the family~{\upshape\ref{comp_diag}} are bounded and
  their omega-limit sets are copies of the base.
\end{theorem}
\begin{proof}
  Proposition~\ref{G1-G2_imply_C1-C4} guarantees that conditions (C1)--(C4)
  are satisfied. For each $i\in\{1,\ldots,m\}$ such that $c_i\equiv 0$, let
  $a_i=-\sup_{\w\in\Om}L_i^+(\w)-1$. Let $A$ be the $m\times m$ diagonal
  matrix with diagonal elements $a_1,\ldots,a_m$. We consider the order
  $\leq_{D,A}$ associated to $\leq_{A,\rho}$ for
  $\rho=\max\{\rho_{11},\ldots,\rho_{mm}\}$. Let us check that the family
  of equations~{\upshape\ref{comp_diag}} satisfies conditions
  {\upshape(F4)} and {\upshape(F6)}. First, let us focus on condition
  (F4). If $c_i\equiv 0$, then $[AD(\w,z)]_i=a_iz_i(0)\leq
  -L_i^+(\w)z_i(0)\leq F_i(\w,y)-F_i(\w,x)$ and we are done. Let us suppose
  that $c_i\not\equiv 0$. Fix $i\in\{1,\ldots,m\}$. Let
  $(\w,x),\,(\w,y)\in\Om\times BU$ with $(\w,x)\leq_{D,A}(\w,y)$ and denote
  $z=y-x$. Then we have that $D_i(\w,z)\geq
  e^{a_i\alpha_i}D_i(\w{\cdot}(-\alpha_i),z_{-\alpha_i})$, whence
  \begin{equation}\label{exp_ord_rho=2alpha}
    z_i(0)-c_i(\w)z_i(-\alpha_i)\geq
    e^{a_i\alpha_i}(z_i(-\alpha_i)-c_i(\w{\cdot}(-\alpha_i))z_i(-\rho_{ii})).
  \end{equation}
  From Theorem~\ref{explicit_D-1}, it follows that $z\geq 0$. Thus
  \[
  \begin{split}
    F_i(\w,y)&-F_i(\w,x)\geq
    -L_i^+(\w)z_i(0)+l_{ii}^-(\w{\cdot}(-\rho_{ii}))z_i(-\rho_{ii})\\
    &+\sum_{j\neq
      i}l_{ij}^-(\w{\cdot}(-\rho_{ij}))z_j(-\rho_{ij})\geq-L_i^+(\w)z_i(0)+l_{ii}^-(\w{\cdot}(-\rho_{ii}))z_i(-\rho_{ii}).
  \end{split}
  \]
  Note that (G3.1) implies that $-a_i-L_i^+(\w)\geq 0$ for all
  $\w\in\Om$. Then, from~\eqref{exp_ord_rho=2alpha}, it follows that
  \[
  \begin{split}
    F_i(\w,y)&-F_i(\w,x)-[AD(\w,z)]_i\geq(-L_i^+(\w)-a_i)z_i(0)+a_ic_i(w)z_i(-\alpha_i)\\
    &+l_{ii}^-(\w{\cdot}(-\rho_{ii}))z_i(-\rho_{ii})\geq[(-L_i^+(\w)-a_i)e^{a_i\alpha_i}-L^+_i(\w)c_i(\w)]z_i(-\alpha_i)\\
    &+[l_{ii}^-(\w{\cdot}(-\rho_{ii}))+(L_i^+(\w)+a_i)c_i(\w{\cdot}(-\alpha_i))e^{a_i\alpha_i}]z_i(-\rho_{ii}).
  \end{split}
  \]
  On the other hand, $D_i(\w{\cdot}(-\alpha_i),z_{-\alpha_i})\geq 0$, which
  in turn implies that
  $z_i(-\alpha_i)-c_i(\w{\cdot}(-\alpha_i))z_i(-\rho_{ii})\geq
  0$. Consequently, thanks to (G3.1) and (G3.2),
  \begin{equation}\label{ineq_trans_rho=2alpha}
    \begin{split}
      F_i(\w,y)&-F_i(\w,x)-[AD(\w,z)]_i\geq\\
      \geq&[l_{ii}^-(\w{\cdot}(-\rho_{ii}))-L_i^+(\w)\,c_i(\w)\,c_i(\w{\cdot}(-\alpha_i))]z_i(-\rho_{ii})\geq
      0.
    \end{split}
  \end{equation}
  As a result, (G3.2) is a sufficient condition for (F4) to hold.\par
  As for (F6), if we had that $D_i(\w{\cdot}s,x_s)<D_i(\w{\cdot}s,y_s)$ for
  all $s\leq 0$, then $z_i(s)>c_i(\w{\cdot}s)z_i(s-\alpha_i)\geq 0$, $s\leq
  0$. It is clear that, if condition (G3.1) is strict, then the first
  inequality in~\eqref{ineq_trans_rho=2alpha} is strict; on the other hand,
  if condition (G3.2) is strict, then the second inequality
  in~\eqref{ineq_trans_rho=2alpha} is strict. This way, the fact that at
  least one of the inequalities in (G3) is strict together with the choice
  of $a_i$ when $c_i\equiv 0$ and an argument similar to the foregoing one
  yield (F6), as expected. The remainder of the theorem follows from
  Theorem~\ref{copy_base_comp}.
\end{proof}
Let us focus on a different approach to the study of
family~\ref{comp_diag}; we will give another valid monotonicity condition
for the transformed exponential order defined on the whole interval
$(-\infty,0]$. For each $i\in\{1,\ldots,m\}$, each $\w\in\Om$ and each
$a\in(-\infty,0]$, we consider $q_{i0}(\w,a)=-L_i^+(\w)-a$ and
\[
\begin{split}
  p_{in}(\w,a)=&-L_i^+(\w)\,c_i^{[n]}(\w)+e^{a(\alpha_i-\rho_{ii})}
  l_{ii}^-(\w{\cdot}(-\rho_{ii}))\,c_i^{[n-1]}(\w{\cdot}(-\rho_{ii})),\\
  q_{in}(\w,a)=&\;q_{in-1}(\w,a)\,e^{a\alpha_i}+p_{in}(\w,a),\,n\in\N.
\end{split}
\]
In the conditions of the next statement, an infinite sequence of
consecutive inequalities is required. If $c$ is Lipschitz continuous, then
the conclusions hold when dealing with Lipschitz continuous initial data.
\begin{theorem}\label{accumulative_coeff}
  Assume hypotheses {\upshape(G1)--(G2)} together with the following one:
  \begin{itemize}
  \item[(G4)] for each $i\in\{1,\ldots,m\}$, if $c_{ii}\not\equiv 0$, then
    $\rho_{ii}\leq\alpha_i$ and there exists $a_i\in(-\infty,0]$ such that,
    for all $\w\in\Om$, there is an $n_0(\w)\in\N\cup\{0\}$ such that
    \[
    \begin{split}
      &q_{in}(\w,a_i)\geq 0,\,n\in\{0,\ldots,n_0(\w)-1\}\text{ (if
      }n_0(\w)\geq 1\text{)},\\
      &q_{in_0(\w)}(\w,a_i)> 0,\,\text{ and}\\
      &p_{in}(\w,a_i)\geq 0\,\text{ for all }n>n_0(\w).
    \end{split}
    \]
  \end{itemize}
  Then, for all $(\w,x)\in\Om\times BU$ such that $\wh D_2(\w,x)$ is
  Lipschitz continuous, the trajectory $\{\tau(t,\w,x):t\geq 0\}$ is
  bounded and its omega-limit set is a copy of the base.
\end{theorem}
\begin{proof}
  First, Proposition~\ref{G1-G2_imply_C1-C4} yields conditions (C1)--(C4).
  For each $i\in\{1,\ldots,m\}$ such that $c_i\equiv 0$, let
  $a_i=-\sup_{\w\in\Om}L_i^+(\w)-1$. Let $A$ be the $m\times m$ diagonal
  matrix with diagonal elements $a_1,\ldots,a_m$ and consider the order
  $\leq_{D,A}$ associated to $\leq_{A,\infty}$. Let us check that the
  family of equations~{\upshape\ref{comp_diag}} satisfies conditions
  {\upshape(F4)} and {\upshape(F6)}. Let $(\w,x),\,(\w,y)\in\Om\times BU$
  with $(\w,x)\leq_{D,A}(\w,y)$ and let $z=y-x$, $\wh z=\wh D_2(\w,z)$. Fix
  $i\in\{1,\ldots,m\}$. As we know
  \[
  z_i(s)=\sum_{n=0}^\infty c_i^{[n]}(\w{\cdot}s)\wh z_i(s-n\alpha_i).
  \]
  Now, if $c_i\equiv 0$, then
  $F_i(\w,y)-F_i(\w,x)-[AD(\w,z)]_i\geq(-L_i^+(\w)-a_i)\wh z_i(0)$ and (F4)
  and (F6) hold. Assume now that $c_i\not\equiv 0$; then we have that
  \[
  \begin{split}
    F_i(\w,y)&-F_i(\w,x)-[AD(\w,z)]_i\geq-L_i^+(\w)z_i(0)+l_{ii}^-(\w{\cdot}(-\rho_{ii}))z_i(-\rho_{ii})-a_i\wh
    z_i(0)\\
    =&-L_i^+(\w)\sum_{n=0}^\infty c_i^{[n]}(\w)\wh z_i(-n\alpha_i)\\
    &+l_{ii}^-(\w{\cdot}(-\rho_{ii}))\sum_{n=0}^\infty
    c_i^{[n]}(\w{\cdot}(-\rho_{ii}))\wh z_i(-\rho_{ii}-n\alpha_i)-a_i\wh
    z_i(0).
  \end{split}
  \]
  Note that $\wh z\geq_A 0$ and $\rho_{ii}\leq\alpha_i$; hence, $\wh
  z_i(-\rho_{ii}-n\alpha_i)\geq e^{a_i(\alpha_i-\rho_{ii})}\wh
  z_i(-(n+1)\alpha_i)$ and
  \begin{equation}\label{eq_monotonicity}
    \begin{split}
      F_i(\w,y)&-F_i(\w,x)-[AD(\w,z)]_i\geq\\
      \geq&(-L_i^+(\w)-a_i)\wh z_i(0)-L_i^+(\w)\sum_{n=1}^\infty c_i^{[n]}(\w)\wh z_i(-n\alpha_i)\\
      &+e^{a_i(\alpha_i-\rho_{ii})}l_{ii}^-(\w{\cdot}(-\rho_{ii}))\sum_{n=1}^\infty
      c_i^{[n-1]}(\w{\cdot}(-\rho_{ii}))\wh z_i(-n\alpha_i)\\
      =&\;q_{i0}(\w,a_i)\wh z_i(0)+\sum_{n=1}^\infty p_{in}(\w,a_i)\wh
      z_i(-n\alpha_i).
    \end{split}
  \end{equation}
  If $n_0(\w)=0$, we are done. Otherwise, we observe that, for all
  $n\in\N\cup\{0\}$, $\wh z_i(-n\alpha_i)\geq e^{a_i\alpha_i}\wh
  z_i(-(n+1)\alpha_i)$. Hence,
  \[
  \begin{split}
    F_i(\w,y)&-F_i(\w,x)-[AD(\w,z)]_i\geq\\
    \geq&\,(q_{i0}(\w,a_i)e^{a_i\alpha_i}+p_{i1}(\w,a_i))\wh
    z_i(-\alpha_i)+\sum_{n=2}^\infty p_{in}(\w,a_i)\wh
    z_i(-n\alpha_i)\\
    =&\,q_{i1}(\w,a_i)\wh z_i(-\alpha_i)+\sum_{n=2}^\infty
    p_{in}(\w,a_i)\wh
    z_i(-n\alpha_i)\\
    \geq&\cdots\geq\,q_{in_0(\w)}(\w,a_i)\wh
    z_i(-n_0(\w)\alpha_i)+\sum_{n=n_0(\w)+1}^\infty p_{in}(\w,a_i)\wh
    z_i(-n\alpha_i).\\
  \end{split}
  \]
  This way, (F4) and (F6) follow easily from (G4). The desired result is a
  consequence of Theorem~\ref{copy_base_comp}.
\end{proof}
The scalar version of the family of functional differential
equations~\ref{comp_diag} with autonomous linear operator and $\rho=\alpha$
was studied in~\cite{arbo} using the standard ordering. In particular, in
this paper it was established that all the minimal sets are copies of the
base. Our next result obtains the same conclusion when the linear operator
is non-autonomous and condition (G5) holds, by using the transformed
exponential ordering.
\begin{proposition}~\label{accumulative_rho=alpha} Assume conditions
  {\upshape(G1)} and {\upshape(G2)}. Consider the
  family~{\upshape\ref{comp_diag}} and assume that, for all
  $i\in\{1,\ldots,m\}$, if $c_i\not\equiv 0$, then $\alpha_i=\rho_{ii}$ and
  the following assertion holds:
  \begin{itemize}
  \item[(G5)] $l_{ii}^-(\w{\cdot}(-\rho_{ii}))-L_i^+(\w)c_i(\w)\geq 0$.
  \end{itemize}
  Then, for all $(\w,x)\in\Om\times BU$, the trajectory
  $\{\tau(t,\w,x):t\geq 0\}$ is bounded and its omega-limit set is a copy
  of the base.
\end{proposition}
\begin{proof}
  Let $\rho>0$ and take the order $\leq_{D,A}$ associated to
  $\leq_{A,\rho}$.  For each $i\in\{1,\ldots,m\}$, let
  $a_i=-\sup_{\w\in\Om}L_i^+(\w)-1$. Let $n_0(\w)=0$; then
  $q_{i0}(\w,a_i)>0$ for all $\w\in\Om$ and
  $p_{in}(\w,a)=c_i^{[n-1]}(\w{\cdot}(-\alpha_i))(l_{ii}^-(\w{\cdot}(-\rho_{ii}))-L_i^+(\w)c_i(\w))\geq
  0$ for all $a\leq 0$ and all $n\in\N$. Following the arguments of
  Theorem~\ref{accumulative_coeff}, we obtain again the
  relation~\eqref{eq_monotonicity}. Now, it is immediate from (G5) that
  condition (G4) holds, which implies the monotonicity of the semiflow. The
  remainder of the proof follows from Theorem~\ref{copy_base_comp}.
\end{proof}
Let us consider the following hypothesis:
\begin{itemize}
\item[(G6)] $\alpha_i>0$, $\rho_{ij}\geq 0$, $0\leq c_i(\w)$ for all
  $i,j\in\{1,\ldots,m\}$ and $\sum_{i=1}^mc_i(\w)<1$ for all $\w\in\Om$.
\end{itemize}
We will give an alternative condition which provides the monotonicity for
the direct exponential order of the semiflow $\tau$ associated to the
family~\ref{comp_diag}. The following result extends the conclusions
in~\cite{krwu} for the scalar periodic case to the $m$-dimensional system
of recurrent NFDEs~\ref{comp_diag}. Note that we provide precise conditions
which assure the monotonicity of the semiflow on $\Om\times BU$. It is
important to mention that, in the present situation, the conclusions
in~\cite{noov} remain valid.
\begin{theorem}\label{like_krwu}
  Assume that conditions {\upshape(G1)} and {\upshape(G6)} hold and,
  moreover, the following condition is satisfied:
  \begin{itemize}
  \item[(G7)] $c$ is continuously differentiable along the flow; let
    $\gamma:\Om\to\R^m$ be its derivative. Besides, for each
    $i\in\{1,\ldots,m\}$, there exists $a_i\in(-\infty,0]$ such that, if
    $A$ is the $m\times m$ diagonal matrix with diagonal elements
    $a_1,\ldots,a_m$ and we consider the order $\leq_A=\leq_{A,\infty}$,
    then, for all $\w\in\Om$, the following inequalities hold:
    \begin{itemize}
    \item[(G7.1)] if $(\w,x),\,(\w,y)\in\Om\times BU$ and $x\leq_Ay$, then
      $F_i(\w,y)-F_i(\w,x)-a_iD_i(\w,y-x)+\gamma_i(\w)(y_i(-\alpha_i)-x_i(-\alpha_i))\geq
      0$.
    \item[(G7.2)] if $(\w,x),\,(\w,y)\in\Om\times BU$ admit a backward
      orbit extension, $x\leq_Ay$ and there exists $J\subset\{1,\ldots,m\}$
      such that $x_i=y_i$ for all $i\not\in J$ and $x_i(s)<y_i(s)$ for all
      $i\in J$ and all $s\leq 0$, then
      $F_i(\w,y)-F_i(\w,x)-a_iD_i(\w,y-x)+\gamma_i(\w)(y_i(-\alpha_i)-x_i(-\alpha_i))>0$
      for all $i\in J$.
    \end{itemize}
  \end{itemize}
  Fix $(\w,x),\,(\w,y)\in\Om\times BU$ such that $x\leq_Ay$. Then
  \[
  u(t,\w,x)\leq_Au(t,\w,y)
  \]
  for all $t\geq 0$ where they are defined. Moreover, all Lipschitz
  continuous initial data give rise to trajectories of the
  family~{\upshape\ref{comp_diag}} which are bounded, and their omega-limit
  sets are copies of the base.
\end{theorem}
\begin{proof}
  Proposition~\ref{G1-G2_imply_C1-C4} guarantees that conditions (C1)--(C4)
  are satisfied. Fix $(\w,x)$, $(\w,y)\in\Om\times BU$ such that
  $x\leq_Ay$. Let $\rho>0$ such that $u(t,\w,x)$, $u(t,\w,y)$ are defined
  on $[0,\rho]$. Let $\eps>0$ and denote by $y^\eps$ the solution of
  \[
  \begin{cases}
    \frac{d}{dt}D(\w{\cdot}t,z_t)=F(\w{\cdot}t,z_t)+\eps{\mathbf 1}, &t\geq
    0,\\
    z_0=y,
  \end{cases}
  \]
  where, ${\mathbf 1}=(1,\ldots,1)^T\in\R^m$. There exists $\eps_0>0$ such
  that, for all $\eps\in[0,\eps_0)$, $z=z({\cdot},\w,x)$ and $y^\eps$ are
  defined on $[0,\rho]$. Let $z^\eps=y^\eps-z$ and denote by $t_1$ the
  greatest element of $[0,\rho]$ such that $z^\eps_{t_1}\geq_A0$. Suppose
  that $t_1<\rho$. Since $z^\eps_{t_1}\geq_A0$, $z^\eps_i(t_1)\geq
  e^{a_i\alpha_i}z^\eps_i(t_1-\alpha_i)$ for all $i\in\{1,\ldots,m\}$ and,
  from (G7.1), it follows that
  \[
  \begin{split}
    \frac{d}{dt}(z^\eps_i(t)&\left.-c_i(\w{\cdot}t)\,z^\eps_i(t-\alpha_i))\right|_{t=t_1}-a_i(z^\eps_i(t_1)-c_i(\w{\cdot}t_1)\,z^\eps_i(t_1-\alpha_i))\\
    &+\gamma_i(\w{\cdot}t_1)\,z^\eps_i(t_1-\alpha_i)=F_i(\w{\cdot}t_1,y^\eps_{t_1})-F_i(\w{\cdot}t_1,z_{t_1})\\
    &-a_i(z^\eps_i(t_1)-c_i(\w{\cdot}t_1)\,z^\eps_i(t_1-\alpha_i))+\gamma_i(\w{\cdot}t_1)\,z^\eps_i(t_1-\alpha_i)+\eps\geq\eps.
  \end{split}
  \]
  Hence, taking $\alpha=\min\{\alpha_1,\ldots,\alpha_m\}$, there exists
  $h\in(0,\alpha)$ such that, if $i\in\{1,\ldots,m\}$ and
  $t\in[t_1,t_1+h]$, then
  \[
  \frac{d}{dt}(z^\eps_i(t)-c_i(\w{\cdot}t)z^\eps_i(t-\alpha_i))-a_i(z^\eps_i(t)-c_i(\w{\cdot}t)z^\eps_i(t-\alpha_i))+\gamma_i(\w{\cdot}t)z^\eps_i(t-\alpha_i)\geq
  0.
  \]
  Let us fix $i\in\{1,\ldots,m\}$; for $t_1\leq s\leq t\leq t_1+h$,
  integrating between $s$ and $t$, we have
  \begin{equation}\label{krwu_ineq}
    \begin{split}
      z^\eps_i(t)-c_i(\w{\cdot}t)z^\eps_i(t-\alpha_i)\geq&
      e^{a_i(t-s)}(z^\eps_i(s)-c_i(\w{\cdot}s)z^\eps_i(s-\alpha_i))\\
      &-\int_s^t\gamma_i(\w{\cdot}u)e^{a_i(t-u)}z^\eps_i(u-\alpha_i)\,du.
    \end{split}
  \end{equation}
  Fix $\eta>0$; there exists an analytic function $\wt
  c_i:[t_1,t_1+h]\to\R$ (for instance, a polynomial) such that
  \[
  0\leq\wt c_i(t)-c_i(\w{\cdot}t)\leq\eta\,\text{ and }\, |\wt
  c_i\,'(t)-\gamma_i(\w{\cdot}t)|\leq\eta
  \]
  for every $t\in[t_1,t_1+h]$.  As a result, from~\eqref{krwu_ineq}, it
  follows that, for $t_1\leq s\leq t\leq t_1+h$,
  \begin{equation}\label{krwu_ineq_analytic}
    \begin{split}
      z^\eps_i(t)&-\wt c_i(t)\,z^\eps_i(t-\alpha_i)\geq
      e^{a_i(t-s)}(z^\eps_i(s)-\wt
      c_i(s)\,z^\eps_i(s-\alpha_i))\\
      &-\int_s^t\wt
      c_i\,'(u)e^{a_i(t-u)}z^\eps_i(u-\alpha_i)\,du-\eta\,\n{z^\eps_{t_1}}_\infty(1+e^{a_i(t-s)}+(t-s)).
    \end{split}
  \end{equation}
  Now, as $\wt c_i$ is analytic, there exist $s=s_0\leq s_1\leq
  s_2\leq\cdots\leq s_J=t$ such that, for all $j\in\{1,\ldots,J\}$, either
  $\wt c_i\,'(u)\geq 0$ or $\wt c_i\,'(u)<0$ for all
  $u\in(s_{j-1},s_j)$. Let $j\in\{0,\ldots,J-1\}$ and fix $N\in\N$ such
  that $N\geq 3$; we define $s_{j0}^N=s_j$,
  $s_{j1}^N=(1-1/N)s_j+1/Ns_{j+1}$, $s_{j2}^N=1/Ns_j+(1-1/N)s_{j+1}$ and
  $s_{j3}^N=s_{j+1}$. Then we have
  \[
  \begin{split}
    \int_{s_j}^{s_{j+1}}&\wt
    c_i\,'(u)e^{a_i(t-u)}z^\eps_i(u-\alpha_i)\,du=\sum_{l=1}^3\int_{s^N_{j\,l-1}}^{s_{jl}^N}\wt
    c_i\,'(u)e^{a_i(t-u)}z^\eps_i(u-\alpha_i)\,du\\
    =&\sum_{l=1}^{3}(\wt c_i(s_{jl}^N)-\wt
    c_i(s^N_{j\,l-1}))e^{a_i(t-u_{jl}^N)}z^\eps_i(u_{jl}^N-\alpha_i)\\
    =&\,\wt c_i(s_{j+1})e^{a_i(t-u_{j3}^N)}z^\eps_i(u_{j3}^N-\alpha_i)-\wt
    c_i(s_{j})e^{a_i(t-u_{j1}^N)}z^\eps_i(u_{j1}^N-\alpha_i)\\
    &+\sum_{l=1}^2\wt
    c_i(s_{jl}^N)\left(e^{a_i(t-u_{jl}^N)}z^\eps_i(u_{jl}^N-\alpha_i)-e^{a_i(t-u_{j\,l+1}^N)}z^\eps_i(u_{j\,l+1}^N-\alpha_i)\right)\\
    \leq&\,\wt
    c_i(s_{j+1})e^{a_i(t-u_{j3}^N)}z^\eps_i(u_{j3}^N-\alpha_i)-\wt
    c_i(s_{j})e^{a_i(t-u_{j1}^N)}z^\eps_i(u_{j1}^N-\alpha_i)
  \end{split}
  \]
  where the points $u_{jl}^N\in[s_{j\,l-1}^N,s_{jl}^N]$ for
  $l=1,\,2,\,3$. As a consequence, $u_{jl}^N-\alpha_i\leq t_1$ for
  $l=1,\,2,\,3$.

  Taking limits when $N\to\infty$, we obtain
  \[
  \begin{split}
    \int_{s_j}^{s_{j+1}}\wt
    c_i\,'(u)e^{a_i(t-u)}z^\eps_i(u-\alpha_i)du\leq&\,\wt
    c_i(s_{j+1})e^{a_i(t-s_{j+1})}\,z^\eps_i(s_{j+1}-\alpha_i)\\
    &-\wt c_i(s_j)e^{a_i(t-s_j)}z^\eps_i(s_j-\alpha_i).
  \end{split}
  \]
  Consequently,
  \[
  \begin{split}
    \int_s^t\wt
    c_i\,'(u)e^{a_i(t-u)}z^\eps_i(u-\alpha_i)\,du\leq&\sum_{j=0}^{J-1}\int_{s_j}^{s_{j+1}}\wt
    c_i\,'(u)e^{a_i(t-u)}z^\eps_i(u-\alpha_i)\,du\\
    \leq&\,\wt c_i(t)z^\eps_i(t-\alpha_i)-\wt
    c_i(s)e^{a_i(t-s)}z^\eps_i(s-\alpha_i)
  \end{split}
  \]
  and, using~\eqref{krwu_ineq_analytic}, it yields
  \[
  z^\eps_i(t)-e^{a_i(t-s)}z^\eps_i(s)\geq-\eta\,\n{z^\eps_{t_1}}_\infty(1+e^{a_i(t-s)}+(t-s)).
  \]
  Letting $\eta\to 0$, we obtain $z^\eps_i(t)-e^{a_i(t-s)}z^\eps_i(s)\geq
  0$ for all $i\in\{1,\ldots,m\}$. This is not possible due to the choice
  of $t_1$. Thus $z^\eps_t\geq_A0$ for all $t\in[0,\rho]$ and, taking
  limits as $\eps\to 0$, $u(t,\w,x)\leq_Au(t,\w,y)$ for all $t\in[0,\rho]$
  and hence for all $t\geq 0$ where they are defined.\par
  Finally, a generalization of the results in~\cite{noov} provides the
  1-covering property under these conditions. The theorem is proved.
\end{proof}
Note that this theorem requires hypothesis (G6), which is significantly
stronger than (G2), which was required for the transformed exponential
order. Besides, according to the previous theory, we should remark that the
application of the direct exponential order requires the differentiability
along the flow of the vector function $c$, instead of just the continuity
of this map, only needed by the transformed exponential order. Thus the
transformed exponential order becomes more natural in the study of NFDEs
with non-autonomous linear $D$-operator.\par
Even in the periodic case, it is known that, given an open set $U\subset
C(\Om,\R^m)$, the subset of differentiable functions is dense, has empty
interior and
\[
\sup\{\n{\mathfrak{c}'}_\infty:\mathfrak{c}\in U\text{ and
}\mathfrak{c}\text{ is differentiable along the flow}\}=\infty
\]
(see Schwartzman~\cite{schw}). As a consequence, the transformed
exponential order is also more advantageous when dealing with rapidly
oscillating differentiable coefficients $c_i$, $i\in\{1,\ldots,m\}$. In
practice, the conditions which allow us to apply the direct exponential
order or the transformed exponential order are frequently quite different;
the particular problem to be studied, will determine the advantages and
disadvantages of each order.\par
We clarify the hypotheses of Theorem~\ref{like_krwu} in some specific
situations. For each $t\in\R$, let $\mathfrak{n}(t)=\min\{t,0\}$. Next, we
give a condition on the coefficients of the equation implying condition
(G7).
\begin{proposition}\label{krwu_specific_conditions}
  Assume that conditions {\upshape(G1)} and {\upshape(G6)} together with
  the following one are satisfied:
  \begin{itemize}
  \item[(G8)] $c$ is continuously differentiable along the flow; let
    $\gamma:\Om\to\R^m$ be its derivative. Besides, for each
    $i\in\{1,\ldots,m\}$, there exists $a_i\in(-\infty,0]$ such that, for
    all $\w\in\Om$, the following inequality holds:
    \[
    -L_i^+(\w)-a_i+\mathfrak{n}(a_ic_i(\w)+\gamma_i(\w))e^{-a_i\alpha_i}>0.
    \]
  \end{itemize}
  Then, for all $(\w,x)\in\Om\times BU$ such that $x$ is Lipschitz
  continuous, the trajectory of $(\w,x)$ for the
  family~{\upshape\ref{comp_diag}} is bounded and its omega-limit set is a
  copy of the base.
\end{proposition}
\begin{proof}
  Let $A$ be the $m\times m$ diagonal matrix with diagonal elements
  $a_1,\ldots,a_m$ and consider the order $\leq_A=\leq_{A,\infty}$. It is
  clear that, if $(\w,x),\,(\w,y)\in\Om\times BU$ and $x\leq_Ay$, then
  \[
  \begin{split}
    F_i(\w,y)&-F_i(\w,x)-a_i(y_i(0)-x_i(0)-c_i(\w)(y_i(-\alpha_i)-x_i(-\alpha_i)))\\
    &+\gamma_i(\w)(y_i(-\alpha_i)-x_i(-\alpha_i))\geq\\
    \geq&(-L_i^+(\w)-a_i+\mathfrak{n}(a_ic_i(\w)+\gamma_i(\w))e^{-a_i\alpha_i})(y_i(0)-x_i(0)).
  \end{split}
  \]
  Consequently, property (G8) guarantees both (G7.1) and (G7.2) and
  Theorem~\ref{like_krwu} yields the expected result.
\end{proof}
We take $\rho_{ii}=2\alpha_i$ under the assumptions of
Proposition~\ref{krwu_specific_conditions}. Note that now condition (G3.2)
is not required. However, if it holds, then (G3.1) is less restrictive than
(G8) even in their autonomous versions.\par
Finally, we turn to the study of~\ref{comp_diag} when another monotonicity
condition is considered, which improves the conclusions of the previous
statement.
\begin{proposition}\label{direct_rho_leq_alpha}
  Assume hypotheses {\upshape(G1)} and {\upshape(G6)} together with the
  following one:
  \begin{itemize}
  \item[(G9)] $c$ is continuously differentiable along the flow; let
    $\gamma:\Om\to\R^m$ be its derivative. Besides, for each
    $i\in\{1,\ldots,m\}$, if $c_i\not\equiv 0$, then
    $\rho_{ii}\leq\alpha_i$ and there exists $a_i\in(-\infty,0]$ such that,
    for all $\w\in\Om$,
    \begin{itemize}
    \item[(G9.1)] $-a_i-L_i^+(\w)\geq 0$,
    \item[(G9.2)]
      $e^{a_i\rho_{ii}}(-a_i-L_{ii}^+(\w))+l_{ii}^-(\w{\cdot}(-\rho_{ii}))
      +e^{a_i(\rho_{ii}-\alpha_i)}\mathfrak{n}(a_ic_i(\w)+\gamma_i(\w))\geq
      0$,
    \end{itemize}
  \end{itemize}
  where at least one of the inequalities is strict. Then all Lipschitz
  continuous initial data for the family~{\upshape\ref{comp_diag}} give
  rise to bounded trajectories and their omega-limit sets are copies of the
  base.
\end{proposition}
\begin{proof}
  First, Proposition~\ref{G1-G2_imply_C1-C4} yields conditions (C1)--(C4).
  For each $i\in\{1,\ldots,m\}$ such that $c_i\equiv 0$, let
  $a_i=-\sup_{\w\in\Om}L_i^+(\w)-1$. Now, let $A$ be the $m\times m$
  diagonal matrix with diagonal elements $a_1,\ldots,a_m$ and consider the
  order $\leq_A=\leq_{A,\infty}$. Let us check that the family of
  equations~{\upshape\ref{comp_diag}} satisfies conditions {\upshape(F4)}
  and {\upshape(F6)}. Let $(\w,x),\,(\w,y)\in\Om\times BU$ with $x\leq_Ay$
  and let $z=y-x$. Fix $i\in\{1,\ldots,m\}$. If $c_i\equiv 0$, then
  $F_i(\w,y)-F_i(\w,x)-[AD(\w,z)]_i+\gamma_i(\w)z_i(-\alpha_i)\geq(-L_i^+(\w)-a_i)z_i(0)$,
  whence (F4) holds and an argument similar to the one given in
  Proposition~\ref{trans_rho=2alpha} yields (F6). Let us assume that
  $c_i\not\equiv 0$; in this case,
  \[
  \begin{split}
    F_i(\w,y)&-F_i(\w,x)-[AD(\w,z)]_i+\gamma_i(\w)z_i(-\alpha_i)\geq-L_i^+(\w)z_i(0)\\
    &+l_{ii}^-(\w{\cdot}(-\rho_{ii}))z_i(\w{\cdot}(-\rho_{ii}))-a_i(z_i(0)-c_i(\w)z_i(-\alpha_i))+\gamma_i(\w)z_i(-\alpha_i)\\
    =&(-a_i-L_i^+(\w))z_i(0)+l_{ii}^-(\w{\cdot}(-\rho_{ii}))z_i(-\rho_{ii})+\mathfrak{n}(a_ic_i(\w)+\gamma_i(\w))z_i(-\alpha_i).
  \end{split}
  \]
  Since $z\geq_A 0$ and $\rho_{ii}\leq\alpha_i$, we have that
  $z_i(-\rho_{ii})\geq e^{a_i(\alpha_i-\rho_{ii})}z_i(-\alpha_i)$,
  $z_i(0)\geq e^{a_i\rho_{ii}}z_i(-\rho_{ii})$ and, thanks to (G9.1), it
  follows that
  \begin{equation}\label{ineq_direct_rho_leq_alpha}
    \begin{split}
      F_i(\w,&\,y)-F_i(\w,x)-[AD(\w,z)]_i+\gamma_i(\w)z_i(-\alpha_i)\geq[e^{a_i\rho_{ii}}(-a_i-L_i^+(\w))\\
      &+l_{ii}^-(\w{\cdot}(-\rho_{ii}))+e^{a_i(\rho_{ii}-\alpha_i)}\mathfrak{n}(a_ic_i(\w)+\gamma_i(\w))]z_i(-\rho_{ii}).
    \end{split}
  \end{equation}
  As a result, (G9.2) is a sufficient condition for (F4) to hold.\par
  As for (F6), if we had that $D_i(\w{\cdot}s,x_s)<D_i(\w{\cdot}s,y_s)$ for
  all $s\leq 0$, then $z_i(s)>c_i(\w{\cdot}s)z_i(s-\alpha_i)\geq 0$, $s\leq
  0$.  Clearly, if condition (G9.1) is strict, then the first inequality
  in~\eqref{ineq_direct_rho_leq_alpha} is strict. As a result, the fact
  that at least one of the inequalities in (G9) is strict together with an
  argument similar to the previous one yield (F6), as desired. The rest of
  the proof follows from Theorem~\ref{like_krwu}.
\end{proof}
Observe that, under the hypotheses of
Proposition~\ref{direct_rho_leq_alpha}, only two supplementary conditions
are required, instead of the infinite sequence needed in
Theorem~\ref{accumulative_coeff} to apply the transformed exponential
order.

\end{document}